\documentclass[11pt]{article}
\usepackage{amsmath,amsthm,amsfonts,latexsym,amssymb,supertabular,amscd}
\usepackage{a4}
\usepackage{hyperref}
\usepackage[titletoc,title]{appendix}

\def\Q{{\mathbb Q}}
\def\Z{{\mathbb Z}}

\def\R{{\mathbb R}}
\def\C{{\mathbb{C}}}

\def\O{{\mathcal{O}}}

\def\gcd{\mathrm{gcd}}
\def\al{{\alpha}}
\def\be{{\beta}}
\def\ga{{\gamma}}
\def\Ga{{\Gamma}}
\def\de{{\delta}}
\def\la{{\lambda}} 
\def\La{{\Lambda}}
\def\th{{\theta}}
\def\ep{{\epsilon}}

\def\om{{\omega}}

\def\ida{{\mathfrak a}}

\def\idp{{\mathfrak p}}
\def\idP{{\mathfrak P}}
\def\ordp{\mathrm{ord}_p}

\def\rd{\mathrm{Rad}}
\def\norm{\mathrm{Norm}}
\def\proof{{\bf Proof}.\:}

\newcommand{\proofend}{\hspace*{1mm} \hfill{$\Box$}}
\newcommand{\Nm}[2]{\mathrm{N}_{#1/#2}}
\newcommand{\ideal}[1]{\langle #1\rangle}
\newcommand{\idclass}[1]{[#1]}

\newcommand{\cnj}[2]{#1^{(#2)}}
\newcommand{\ord}[2]{\mathrm{ord}_{#1}(#2)}
\newcommand{\ekth}[2]{{\mathrm{v}_{#1}(#2)}}
\newcommand{\tbsp}{\rule{0pt}{15pt}}

\newcommand{\ndist}[1]{\parallel\!\!#1\!\!\parallel}

\newtheorem{proposition}{Proposition}[section]
\newtheorem{theorem}[proposition]{Theorem}
\newtheorem{lemma}[proposition]{Lemma}

\numberwithin{equation}{section}

\begin{document}

\baselineskip=17pt

\title{Complete solution of the Diophantine Equation $x^{2}+5^{a}\cdot 11^{b}=y^{n} $}
\author{
G.~Soydan \thanks{Department of Mathematics, Uludag University, Bursa, Turkey}
\and
N.~Tzanakis \thanks{Department of Mathematics, University of Crete, Greece}
}
\date{} %% 9 July  2016}

\maketitle

\renewcommand{\thefootnote}{}



\renewcommand{\thefootnote}{\arabic{footnote}}

\setcounter{footnote}{0}
\begin{abstract}
\noindent
The title equation is completely solved in integers $(n,x,y,a,b)$, where $n\geq 3$, $\gcd(x,y)=1$ and
$a,b\geq 0$. The most difficult stage of the resolution is the explicit resolution of a quintic Thue-Mahler 
equation. Since it is for the first time -to the best of our knowledge- that such an equation is solved in
the literature, we make a detailed presentation of the resolution; this gives our paper  also an expository 
character. 
\end{abstract}
 %
% \textup{2010} \textit{Mathematics Subject Classification}: \textup{11D25}
\section{Introduction}
            \label{sec intro}
The title equation belongs to the general class of Diophantine equations of the form
\begin{equation}\label{eq.1}
x^{2}+D=y^{n}, \quad x,y\geq1, \quad n\geq 3,
\end{equation}
where $D$ is a positive integer all whose prime factors belong to a finite set $S$ of at least two distinct
primes.
All solutions of the Diophantine equation \eqref{eq.1} have been determined for various sets $S$: 
In \cite{Lu} for $S=\{2,3\}$, in \cite{LT1} for $S=\{2,5\}$, in \cite{CDLPS} for $S=\{2,11\}$, in 
\cite{LT2} for $S=\{2,13\}$, in \cite{Dabrowski} for $S=\{2,17\}$, $S=\{2,29\}$, $S=\{2,41\}$, in 
\cite{SUZ} for $S=\{2,19\}$.
Note that, in all these cases, $S=\{2,p\}$, where $p$ is an odd prime. The case of $S=\{2,p\}$, with a 
general odd prime $p$, was recently studied by H.~Zhu, M.~Le, G.~Soydan and A.~Togb\'{e} \cite{ZLST}, 
who gave all the solutions of   
$ x^{2}+2^{a}p^{b}=y^{n}, \; x\geq1, y>1,\; \gcd(x,y)=1,\; a\geq0, b>0, \; n\geq3 $
under some conditions.

Several papers deal with the Diophantine equation \eqref{eq.1} when $S$ 
contains at least two distinct odd primes. Thus, all solutions of the Diophantine equation \eqref{eq.1} 
were given in \cite{ALT} for $S=\{5,13\}$, in \cite{PR} for $S=\{5,17\}$, in \cite{S1}, \cite{S2} for 
$S=\{7,11\}$ -except for the case when $ax$ is odd and $b$ is even-, in \cite{BP} for $S=\{11,17\}$, in 
\cite{GLT} for $S=\{2,5,13\}$, in \cite{CDILS} for $S=\{2,3,11\}$, in \cite{GMT} for $S=\{2,5,17\}$. 
In \cite{P}, Pink gave all the non-exceptional solutions of the equation \eqref{eq.1} 
(according to the terminology of that paper) for $S=\{2,3,5,7\}$. A survey of these 
and many others can be found in \cite{BP}, \cite{BM}. Very recently, the equations with $S=\{2,3,17\}$
and $S=\{2,13,17\}$ were solved in \cite{GMT2}.

In \cite{CDST}, I.N.~Cangul, M.~Demirci, G.~Soydan and N.~Tzanakis gave the complete solution  
$(n,a,b,x,y)$ of the Diophantine equation \eqref{eq.1} for $S=\{5,11\}$ when $\gcd(x,y)=1$, except for the 
case when $abx$ is odd. In this paper we treat this remaining case, proving thus the following:
\begin{theorem}\label{thm main theorem}
For the integer solutions of the equation 
\begin{equation} \label{eq title}
x^2+5^a11^b=y^n, \quad n\geq 3,\; x,y\geq 1,\,\gcd(x,y)=1, \; a,b\geq 0, 
\end{equation}
the following hold:
\\
If $n=3$, the only solutions are: 
$ (a,b,x,y) =  (0,1,4,3),\,(0,1,58,15),\,(0,2,2,5)$, \\ $(0,3,9324,443),\,(1,1,3,4),\,
              (1,1,419,56),\,(2,3,968,99),\,(3,1,37,14),\,(5,5,36599,1226)$.
\\
If $n=4$ there are no solutions $(a,b,x,y)$ and, for $n=6$,  the only integer solution is 
$(a,b,x,y)=(1,1,3,2)$.
\\
If $n= 5$ or $n\geq 7$, the equation has no integer solutions $(a,b,x,y)$.
\end{theorem}
The proof of Theorem \ref{thm main theorem} is already accomplished in \cite[Theorem 1]{CDST} for the 
following cases: (A) $n=3,4,6$, and 
(B) $n\geq 5, n\neq 6$ and either (i) $ab$ is odd and $x$ is even, or (ii) at least one of $a,b$ is even.
Therefore, what remains is:
\begin{quote}
\emph{
To prove that, if $ n\geq 5$ is prime, then the equation \emph{(\ref{eq title})} has no solution 
$(a,b,x,y)$ with $abx$ odd.
}
\end{quote}

This paper has two objectives, the first one being displayed above. 
Second is the systematic discussion in Section \ref{sec main section: n=5} of the resolution of the 
quintic Thue-Mahler equation (\ref{eq ThueMahler5}) which,  
along with the three Appendices (see ``Plan of the paper'' below) lends also an expository character to 
the paper, as it presents in detail the application of the method of N.~Tzanakis \& B.M.M.~de Weger 
\cite{TdW} to the explicit resolution of a \emph{quintic} Thue-Mahler equation. 
To the best of our knowledge, in the literature it is the first example of explicit 
resolution of a \emph{quintic} Thue-Mahler equation. 
Indeed, in \cite{TdW}, the worked example is a \emph{cubic} Thue-Mahler equation; those days --almost 25 
years ago-- the available software was not as developed as to support the application of the 
method to a quintic Thue-Mahler equation; even until today, only very few works are published in 
which Thue-Mahler equations are \emph{explicitly} solved and none of them deals with a quintic 
equation; more specifically:  
In \cite{CDST},  I.N.~Cangul, M.~Demirci, G.~Soydan and N.~Tzanakis need to solve --successfully-- a 
quartic Thue-Mahler equation. In \cite{Kim}, Dohyeong Kim   proposes a method different from that of 
\cite{TdW} --with many examples-- for the explicit resolution of \emph{cubic} Thue-Mahler equations, which 
exploits the modularity of elliptic curves over $\Q$.   
M.A.~Bennett and S.R.~Dahmen \cite{BennDahm}, in their study of \emph{generalized superelliptic} 
equations need to consider some special classes of Thue-Mahler equations. 
These are closely related to the so-called \emph{Klein forms}, which are defined as binary forms 
of the following shape: $F(x,y):=F_n(ax+by, cx+dy)$, where 
$\left(\begin{array}{cc} a & b \\ c & d \end{array}\right)\in\mathrm{GL}_2(\overline{\Q})$, 
$n\in\{2,3,4,5\}$ and $F_2(x,y)=xy(x+y)$, $F_3(x,y)=y(x^3+y^3)$,
$F_4(x,y)=xy(x^4+y^4)$, $F_5(x,y)=xy(x^{10}-11x^5y^5-y^{10})$. Thue-Mahler equations whose left-hand
side is a Klein form are considered. In the case of \emph{cubic} Klein forms, Bennett and Dahmen, 
implemented the method of \cite{TdW}; for the purposes of their paper, explicit resolution of higher 
degree Thue-Mahler equations was not necessary. 
\\
Finally, we mention Kyle Hambrook's M.\,Sc.\,thesis \cite{Hamb}, where the method of \cite{TdW} is
revisited, certain improvements are included and, most importantly, a long\footnote{More than 200 pages!}
{\sc magma} program is developed for the automatic resolution of the general Thue-Mahler equation; the 
program needs as its input only the coefficients and the primes of the equation. 
No examples of Thue-Mahler equations of degree greater than three are discussed. 
As we checked, the program runs very succesfully with ``reasonable'' cubic Thue-Mahler equations. 
This work is, certainly, a good contribution to the project of the automatic resolution of Thue-Mahler 
equations.
However, in the case of our quintic Thue-Mahler equation (\ref{eq ThueMahler5}),
it took 72 days on an Apple computer with the following characteristics: 
Processor Intel i5, 2.5 GHz, 4GB RAM, 1600 MHz DDR3.
Therefore, we preferred to develop also our own {\sc magma} program, 
far less automatic than that of Hambrook, which needs human intervention at various points. With this
program, the resolution of the quintic Thue-Mahler equation took less than $2\,\mbox{h}\,20'$.
Besides this huge difference in computation time, this ``primitive'' type of computer-aided resolution, 
however, has the advantage that it allows a rather transparent presentation of the very complicated 
resolution.
It is our belief that the experience in the ``technical details'' needed
for the development of a very satisfactory ``Thue-Mahler automatic solver'' requires the resolution 
of quite a number of specific Thue-Mahler equations (over $\Z$) of degree $\leq 6$ (at least), with 
corresponding number fields of various types.\footnote{We mention here that very satisfactory
automatic Thue solvers are included, for example, in {\sc pari} \cite{pari} and {\sc magma} \cite{magma} 
since long time, and are based on Bilu-Hanrot's improvement \cite{BiHa} of Tzanakis-de Weger 
\cite{TdW-Thue} method for solving Thue equations. The development of an automatic Thue solver is, 
certainly, a difficult job but, anyway, much easier than an analogous job for a Thue-Mahler equation. We 
use the {\sc magma} Thue solver in Subsection \ref{subsubsec final}.}

\emph{Plan of the paper}. In Section \ref{sec n>5} we prove that equation (\ref{eq title}) has no 
solutions with $abx$ odd and prime $n\geq7$. Thus, we are reduced to proving that our equation  has no 
solutions with $abx$ odd and $n=5$. This is accomplished in Section \ref{sec main section: n=5} which
is the heart of the paper and is divided into two subsections.
In Subsection \ref{subsec reduce to a TM eq}, using standard algebraic Number Theory, we reduce  
the equation $x^2+5^a11^b=y^5$ -with $\gcd(x,y)=1$ and $abx$ odd- to the quintic Thue-Mahler 
equation (\ref{eq ThueMahler5}), whose right-hand side is $-2^53^4 5^{z_1}11^{z_2}$, where 
$z_1=(a-1)/2$ and $z_2=(b-1)/2$ are now our non-negative unknown integers. 
Then, Subsection \ref{subsec solution of T-M eq} is devoted to the resolution of that 
Thue-Mahler equation, quite a complicated task. In order to make the exposition of our resolution as clear 
as possible, we divided Subsection \ref{subsec solution of T-M eq} into nine (sub)subsections from 
(sub)Subsection \ref{subsubsec arithm data F} through (sub)Subsection \ref{subsubsec final}. 
\\ \textbullet\,
 In Subsection \ref{subsubsec arithm data F}, using standard arguments from algebraic Number 
Theory along with the valuable routines of {\sc magma} \cite{magma}, we reduce our quintic Thue-Mahler 
equation to the ideal equations (\ref{eq alternative 1}) and (\ref{eq alternative 2}), in the right-hand 
side of which appear the unknown non-negative integers $z_1$ and $z_2$. 
\\ \textbullet\,
In Subsection \ref{subsubsec working 5adically}, working $5$-adically, we prove that $z_1\leq 27$ 
(implying $a\leq 55$). 
\\ \textbullet\,
In Subsection \ref{subsubsec working 11adically} we we work $11$-adically. Making also use of the upper 
bound $z_1\leq 27$, we are led to the following situation:
Instead of solving one (quintic) Thue-Mahler equation in which the exponents of the 
primes 5 and 11 are among the unknowns, we are led to solving 28  similar Thue-Mahler equations, with all 
having the same left-hand side and right-hand sides in which only the exponent of
the prime 11 is among the unknowns. This is certainly a profit; these 28 equations can be treated as one
equation, namely, equation (\ref{eq final Thue-Mahler}). In this last equation, besides the unknown 
integer $z_2=(b-1)/2$, three more unknown integers $a_1,\ldots,a_4$ make their appearance; these are
the exponents of the four fundamental units of the quintic field related to the Thue-Mahler equation.
\\ \textbullet\,
In Subsection \ref{subsubsec S-unit eq} the aforementioned equation (\ref{eq final Thue-Mahler}) leads 
to the three-term $S$-unit equation (\ref{eq (12) of TdW}); this is the basic step towards the 
use of \emph{Linear Forms in Logarithms} in both the real/complex and the $p$-adic sense; consequently:
\\ \textbullet\, 
In Subsection \ref{subsubsec n1<c13(logH+c14)} we are in a position to apply a powerful 
result of Kunrui Yu (Theorem \ref{thm Yu thm simplified} in this paper) which, given the algebraic numbers 
$\al_1,\ldots,\al_n$ and a prime $p$, provides an upper bound for the $p$-adic valuation of 
$\al_1^{b_1}\cdots\al_n^{b_n} - 1$, for any $b_1,\ldots,b_n\in\Z$,  in terms of 
$\log \max\{3, |b_1|,\ldots,|b_n|\}$.\footnote{This is, actually, equivalent to giving a lower bound of the 
$p$-adic absolute value for a linear form in the \emph{$p$-adic logarithms} of $\al_1,\ldots,\al_n$.} 
Combining the result of this application with the instructions in  p.~238 of \cite{TdW} we manage to  
bound $z_2$ in terms of $\log\max\{z_2,|a_1|,\ldots,|a_4|\}$. 
We remark here that, in order to conform with the notation of \cite{TdW}, the recipes of which we follow 
very closely, we denote $z_2$ by $n_1$. 
\\ \textbullet\,
In Subsection \ref{subsubsec 1st ub} we apply another strong result due to 
E.M.~Matveev (Theorem \ref{thm Matv} in this paper) which, given the algebraic numbers 
$\al_1,\ldots,\al_n$, provides\footnote{Under one or two mild conditions.} a lower bound for 
$|b_1\log\al_1+\cdots+b_n\log\al_n|$, for any $b_1,\ldots,b_n\in\Z$, in terms of 
$\log \max\{3, |b_1|,\ldots,|b_n|\}$. Applying this theorem in our case and combining the result with the
detailed instructions of \cite[\textsection\textsection\,10,11]{TdW}, we obtain a numerical upper
bound for $H=\max\{n_1,|a_1|,|a_2|,|a_3|,|a_4|\}$; see (\ref{1st bound H}). Then, this upper bound of $H$,
in combination with the result of the previous Subsection \ref{subsubsec n1<c13(logH+c14)}, gives a
specific numerical upper bound for $n_1$, which is considerably smaller than $H$; 
see (\ref{1st bound n1}). 
However, both upper bounds are huge, of the size of $10^{43}$ and $10^{32}$, respectively and need to
be reduced to a manageable size, as discussed in \cite[Section 13]{TdW}. 
\\ \textbullet\,
In Subsection \ref{subsubsec padic reduc} we apply the so-called ``$p$-adic reduction step'', which is 
described in detail in Sections 12,14 and 15 of \cite{TdW} and reduce the upper bound for $n_1$ to 
$n_1\leq 207$. 
\\ \textbullet\,
In Subsection \ref{subsubsec real reduc} we combine this extremely smaller upper bound with the 
(still remaining) huge upper bound of $H$, and do the ``real reduction step'', following the instructions 
of \cite[Section 16]{TdW}. With this step we get the bound $H\leq 231$. 
\\ \textbullet\,
In the final Subsection \ref{subsubsec final} we discuss how we proceed with a further reduction, by 
successively repeating the ``$p$-adic reduction'' and the ``real reduction'' steps two more times, until 
we obtain the bound $(z_2=)\,n_1\leq 21$ and $H\leq 34$.  At this point we don't need the bound 
$H\leq 34$; as explained in Subsection \ref{subsubsec final}, we are left with the task of solving 560 
\emph{Thue} equations (\ref{eq ThueMahler5}), whose right-hand sides runs through the set  
$\{-2^5 3^4 5^{z_1}11^{z_2} : 0\leq z_1\leq 27,\, z_2=0 \;\mbox{or}\; 3\leq z_2\leq 21\}$; 
note that the left-hand sides of all these Thue equations  are identical. 
For their solution we use {\sc magma}'s implementation of Bilu \& Hanrot's method \cite{BiHa}.  
It turns out that no solutions exist and this completes the proof that equation (\ref{eq title}) with 
$abx$ odd and $n=5$ has no solutions. Since in Section \ref{sec n>5} we have also proved that equation 
(\ref{eq title}) with $abx$ odd and $n>5$ has no solutions, we have completed the proof of 
Theorem \ref{thm main theorem}.

In the Appendices \ref{Append working padically} through  \ref{Append field K} at the end of the paper 
we collect some theoretical facts and give some information about how these are realized in practice 
with the use of {\sc magma} \cite{magma}. We also give the results of a few computations. The huge 
algebraic numbers in Appendix \ref{Append field K} are not strictly necessary; however, they are useful in 
giving the reader a sense of what ``monsters'' are involved in such a task. 
We hope that the appendices will make transparent our way of work and friendly the reading of our paper.    

\textbf{Acknowledgements}.
This work started during the second author's visit to Department of Mathematics of Uluda\u{g} University 
which was supported by Visiting Scientist Fellowship Programme (2221) of the Scientific and 
Technological Research Council of Turkey (T\"{U}B\.{I}TAK). The first author was supported by the 
Research Fund of Uluda\u{g} University under Project No: F-2015/23, F-2016/9, and would like to thank to 
Professor Mike Bennett for useful suggestions.

The authors thank the anonymous referee for his careful reading and useful comments.
\section{Equation (\ref{eq title}) with $abx$ odd and prime $n\geq 7$}
                \label{sec n>5}
\begin{proposition} \label{prop n<7}
Equation \emph{(\ref{eq title})} has no solutions with $xab$ odd and prime $n\geq 7$.
\end{proposition}
\proof
We assume that a solution $(x,a,b,n,y)$ in which $xab$ is odd and $n$ is a
prime $\geq 7$ does exist, and we put $a=na_1+\al$ and $b=nb_1+\be$, where 
$0\leq\al,\be<n$, so that our equation becomes
\begin{equation} \label{eq initial}
5^{\al}11^{\be}(-5^{a_1}11^{b_1})^n +y^n = x^2.
\end{equation}
Without loss of generality we assume that $x\equiv 1\pmod 4$.
According to the notation etc of \cite[Section 14]{Siksek1}, this is a ternary 
equation of signature $(n,n,2)$, so that it falls under the scope of the recipe
described in \cite[\textsection 14.2]{Siksek1}.
Accordingly, we have the following table which shows how we will apply that
recipe in our case.
\\[-5mm]
\begin{center}
\tablecaption{Application of the recipe in  \textsection 14.2 of \cite{Siksek1}} 
 \label{Table apply recipe p,p,2}  
%%%%%%%%%%%%%%%%%%%%%%%%%%%%%%%%%%%%%%%%%%%%%%%%%%%%%%%%%%%%%%%%%%%%%%%%%%%%%%%%%%%%%%%
\tablefirsthead{ 
\multicolumn{2}{c}{} \\[-8mm] \hline
Notations/conditions in \cite[\textsection 14.2]{Siksek1}  
                        & Interpretations in this paper
                                                             \\  \hline\hline
               }
\tablehead{ \hline
\multicolumn{2}{|l|}{{\small \sl continued from previous page}} \\ \hline
Notations/conditions in \cite[\textsection 14.2]{Siksek1}  
                        & Interpretations in this paper
                                                             \\
          }
 \tabletail{\hline \multicolumn{2}{|r|}{\small\sl continued on next page}\\ \hline}
 \tablelasttail{\hline\hline}
 \begin{supertabular}{|c|c|}
   & \\[-10pt]
$A$ & $5^{\al}11^{\be}$ \\ \hline
  & \\[-10pt]
$B$ & $1$ \\ \hline
%  & \\[-10pt]
$C$ & $1$ \\ \hline
  & \\[-10pt]
$x$ & $-5^{a_1}11^{b_1}$ \\ \hline
$y$ & $y$ \\ \hline
$z$ & $x$ \\ \hline
general prime $q$  & general prime $q$ \\ \hline
$p$ & $n$ \\ \hline
 % & \\[-8pt]
$\ord{q}{B}<p$ & trivially satisfied \\ \hline
 % & \\[-8pt]
$\ord{q}{A}<p$ & $\al,\be <p$ \\ \hline
%  & \\[-10pt]
$C$ square-free & trivially satisfied \\
 \end{supertabular}
 \end{center}
Since $y$ is even, $x\equiv 1\pmod 4$ and $n\geq 7$, our equation falls in 
case (v) of \cite[\textsection 14.2]{Siksek1} and we deal with the elliptic 
curve
\[
E_3: Y^2+XY=X^3+\frac{x-1}{4}X^2+\frac{y^n}{64}X.
\]
According to a result of Bennett and Skinner \cite[Lemma 2.1]{BennSkin} 
(or \cite[Theorem 16]{Siksek1}), the discriminant and conductor of this 
elliptic curve are, respectively
\[
\Delta_3=-2^{-12}5^{\al}11^{\be}(5^{a_1}11^{b_1}y^2)^n =
-2^{-12}5^{a}11^{b}y^{2n}, \quad 
N_3=55\,\rd(y)=5\cdot 11\prod_{q|y}\!q,
\]
where in the last product $q$ is prime. 
According to \cite[Theorem 16\,(c)]{Siksek1}, there exists a newform $f$ of
level $N_n=2\cdot 5\cdot 11=110$, such that $E_3\thicksim_n f$ 
(\emph{$E_3$ arises from $f\bmod\,n$}; see \cite[\textsection 5]{Siksek1}).

A computation using {\sc magma} returns three rational newforms of level $110$,
namely\footnote{Below $q$ denotes the ``$q$-variable'' of the modular form and has 
nothing to do with primes.}
\begin{eqnarray*}
f_1 & = &
q - q^2 + q^3 + q^4 - q^5 - q^6 + 5q^7 - q^8 - 2q^9 + q^{10} + q^{11} + O(q^{12}) \\
f_2 & = & 
q + q^2 + q^3 + q^4 - q^5 + q^6 - q^7 + q^8 - 2q^9 - q^{10} - q^{11} + O(q^{12}) \\
f_3 & = & 
  q + q^2 - q^3 + q^4 + q^5 - q^6 + 3q^7 + q^8 - 2q^9 + q^{10} + q^{11} + O(q^{12})        
\end{eqnarray*}
and a non-rational newform
\[
f_4= q - q^2 + \al q^3 + q^4 + q^5 - \al q^6 - \al q^7 - q^8 
+ (5-\al)q^9 - q^{10} - q^{11} + O(q^{12}), 
\]
where $\al^2+\al-8=0$, along with its conjugate newform.

Now we apply \cite[Proposition 9.1]{Siksek1} to $E=E_3$ and $f=f_i$, 
$i=1,2,3,4$. 
Our notation refers to that Proposition. Since $(X,Y)=(0,0)$ is a $2$-torsion point 
on $E_3$, we take $t=2$. Also $N'=110$ and we choose the prime $\ell=3$, noting 
that $\ell\nmid N'$ and $\ell^2\nmid N_3$. Then, 
\[
S_3=\{a\in\Z: -2\sqrt{3}\leq a\leq 2\sqrt{3},\;\mbox{$a$ is even}\}=
\{-2,0,2\}.
\]
Also, denoting by $c_{3i}$ the coefficient of $q^3$ in the newform $f_i$,
we compute
\[
B_3(f_i)=(4^2-c_{3i}^2)\prod_{a\in S_3}(a-c_{3i}) =
\begin{cases}
3^2\cdot 5 & \mbox{if $i=1,2$} \\
-3^2\cdot 5 & \mbox{if $i=3$.}
\end{cases}
\]
For the newform $f_4$ we compute
\[
B_3(f_4)= 3\cdot\mathrm{Norm}_{\Q(\al)/\Q}(4^2-\al^2)
     \prod_{k=-1}^1\mathrm{Norm}_{\Q(\al)/\Q}(2k-\al)=-3^3\cdot 2^9.
\]
According to the conclusion of \cite[Proposition 9.1]{Siksek1}, $n$ must divide
$B_3(f_i)$ for some $i\in\{1,2,3,4\}$, which is impossible since we assumed
that $n$ is a prime $\geq 7$.
\proofend
\section{Equation (\ref{eq title}) with $abx$ odd and $n=5$ }
  \label{sec main section: n=5}
 \subsection{Reduction to the Thue-Mahler equation (\ref{eq ThueMahler5})}  
   \label{subsec reduce to a TM eq}
In view of Proposition \ref{prop n<7}, we are left with $n=3,5$.
The case $n=3$ is already solved completely; see \cite[Proposition 2]{CDST}.
In particular, the only solutions with $abx$ odd are the following:
\[
(a,b,x,y)=(1,1,3,4),\,(1,1,419,56),\,(3,1,37,14),\,(5,5,36599,1226).
\]

It remains to treat the title equation when $n=5$ and $xab$ is odd.
We write our equation
\begin{equation} \label{eq initial small p}
x^2+55z^2=2^5y_1^5,\quad x\equiv 1\hspace{-3mm}\pmod 4,\; y=2y_1,\; z=5^{(a-1)/2}11^{(b-1)/2}
\end{equation}
and work in the field
\[
L=\Q(\rho),\quad \rho =\frac{1+\sqrt{-55}}{2}\quad (\rho^2-\rho+14=0).
\]
Using either {\sc pari-gp} \cite{pari} or {\sc magma} \cite{magma} we can obtain the following facts
about the number-field $L$:
\begin{itemize}
\item The class-number is $4$.
\item $\ideal{2}=\idp_2\idp_2'$, where
       $\idp_2=\ideal{2,\rho}$, $\idp_2'=\ideal{2,3+\rho}$.
\item The order of the ideal-class of both $\idp_2$ and $\idp_2'$ in the 
ideal-class group is $4$. More specifically,
$\idp_2^4=\ideal{2-\rho}$ and $\idp_2'=\ideal{1+\rho}$.   
\item $\ideal{\rho}=\idp_2\ideal{7,\rho}$.   
\end{itemize}
From (\ref{eq initial small p}) we obtain the ideal equation
\begin{equation}  \label{eq initial ideal}
\idp_2^5\idp_2'{^{5}}\ideal{y_1}^5 =\ideal{x+z\sqrt{-55}}\ideal{x-z\sqrt{-55}}=
\ideal{x-z+2z\rho}\ideal{x+z-2z\rho}
\end{equation}
Our first observation is that no prime ideal factor over an odd rational prime
can divide both ideal factors in the right-hand side of (\ref{eq initial ideal}).
Indeed, if $\idp$ is a prime ideal over an odd rational prime, then $\idp|2x$,
hence $\idp|x$. But $\idp|\ideal{x+z\sqrt{-55}}$, hence $\idp |z\sqrt{-55}$.
It follows that $\idp|55$, which contradicts $\gcd(x,55)=1$.
\\
Next we observe that $2$ divides both $x-z+2z\rho$ and $x+z-2z\rho$.
From $\ideal{\rho}=\idp_2\ideal{7,\rho}$ we see that 
$\ideal{2z\rho}=\idp_2^2\idp_2'\times\mbox{(ideal relatively prime to $2$)}$.
Also $x+z\equiv 2\pmod{4}$ shows that $\ord{\idp_2}{x+z}=1=\ord{\idp_2'}{x+z}$.
As a consequence, $\ord{\idp_2}{x+z-2z\rho}=1=\ord{\idp_2'}{x+z-2z\rho}$, which
implies that $\idp_2\|\ideal{x+z-2z\rho}$ and $\idp_2'\|\ideal{x-z+2z\rho}$.
Similarly, starting from $x-z\equiv 0\pmod{4}$, we conclude that
$\idp_2'\|\ideal{x-z-2z\rho}$. 

Combining the above small observations with (\ref{eq initial ideal}), we conclude
that
\begin{equation} \label{eq two ideal factors}
\ideal{x-z+2z\rho}=\idp_2'\idp_2^4\ida_1^5, \quad
\ideal{x+z-2z\rho}=\idp_2\idp_2'{^4}\ida_2^5,
\end{equation}
where $\ida_1,\ida_2$ are relatively prime (integral) ideals, such that
$\ida_1\ida_2=\ideal{y_1}$.

The first equation (\ref{eq two ideal factors}) becomes
\[
\ideal{x-z+2z\rho}=\idp_2'\ida_1^5\ideal{2-\rho}.
\]
Since the class-number is $4$, the above ideal equation implies the ideal-class
equation 
$\idclass{1}=\idclass{\idp_2'}\idclass{\ida_1}=\idclass{\idp_2'\ida_1}$.
Therefore, $\idp_2'\ida_1$ is a principal ideal, so that, on multiplying both 
sides of the above displayed equation by $\idp_2'{^4}$ and setting 
$\idp_2'\ida_1=\ideal{u+v\rho}$ ($u,v\in\Z$) we finally arrive at the following
(element) equation
\begin{equation} \label{eq when n=5}
(1+\rho)(x-z+2z\rho)=(2-\rho)(u+v\rho)^n.
\end{equation}
In (\ref{eq when n=5}) we equate coefficients of $\rho$ in both sides,
as well as rational parts, obtaining thus the following two relations:
\[
3u^5+65vu^4-290v^2u^3-2110v^3u^2+975v^4u+3149v^5 
             =  -32\cdot 5^{(a-1)/2}11^{(b-1)/2},
\]
\[ 
23u^5-355vu^4-3930v^2u^3+6010v^3u^2+30515v^4u-2311v^5 = -32x.             
\]
By multiplying both sides of the first displayed equation by $3^4$ we
get 
\begin{eqnarray*} 
\lefteqn{ \hspace{-20mm}
(3u)^5+65v(3u)^4-870v^2(3u)^3-18990v^3(3u)^2+26325v^4(3u)+255069v^5 
        }
             \nonumber \\
   & &      \hspace*{40mm} = -2^5\cdot 3^4\cdot 5^{(a-1)/2}11^{(b-1)/2}
\end{eqnarray*}  
In order to conform precisely with the notations of \cite{TdW} the  method of which we will apply in this 
section, we set 
\[
(3u,v)=(x,y),\quad z_1=(a-1)/2,\quad z_2=(b-1)/2,
\]
so that
\begin{equation}  \label{eq ThueMahler5}
\norm_{F/\Q}(x-y\th) = -2^5\cdot 3^4\cdot 5^{z_1}11^{z_2},
\end{equation}
where $ F=\Q(\th)$ with $g(\th)=0 $ and 
\begin{equation} \label{eq def g}
g(t)=t^5+65t^4-870t^3-18990t^2+26325t+255069,
\end{equation}
with (polynomial) discriminant $D_{\th}=2^{32}3^{12}5^{11}11^6$.
\subsection{Resolution of the Thue-Mahler equation (\ref{eq ThueMahler5})} 
    \label{subsec solution of T-M eq}
\subsubsection{From equations (\ref{eq ThueMahler5}) to ideal equations (\ref{eq alternative 1}) and (\ref{eq alternative 2})} 
                                                        \label{subsubsec arithm data F}
We need the following arithmetical data for the number field $F$.
\begin{itemize}
\item $F$ is a totally real field with class-number $1$.
\item An integral basis is $1,\be_2,\be_3,\be_4,\be_5 $, where
\[
\be_2=\textstyle{\frac{1}{2}}(\th+1), \quad 
     \be_3= \textstyle{\frac{1}{24}}(\th^2 + 2\th + 9),
\]
\[
\be_4=\textstyle{\frac{1}{7920}}(\th^3 + 227\th^2 + 3603\th + 3969), \quad
\be_5=\textstyle{\frac{1}{95040}}(\th^4 + 8\th^3 + 1410\th^2 + 46512\th + 9909).
\]
\item A quadruple of fundamental units is the following:
\begin{eqnarray*}
\ep_1 & = & 
  \textstyle{\frac{1}{15840}}(\th^4 + 62\th^3 - 852\th^2 - 1806\th + 6435)
\quad (\norm(\ep_1)=1) \\
\ep_2 & = & 
\textstyle{\frac{1}{95040}}(-\th^4 - 104\th^3 + 558\th^2 + 35280\th + 108027)
\quad (\norm(\ep_2)=-1) \\
\ep_3 & = & 
\textstyle{\frac{1}{23760}}(\th^4 + 77\th^3 + 243\th^2 + 99\th - 10260)
\quad (\norm(\ep_3)=1) \\
\ep_4 & = & 
\textstyle{\frac{1}{95040}}(7\th^4 + 596\th^3 + 5730\th^2 - 25596\th - 210897)
\quad (\norm(\ep_4)=1)
\end{eqnarray*}
\item
Prime factorization of $2$: 
\begin{align*}
2 & =-\ep_1^2\ep_2^2\ep_3\ep_4^2\pi_2^5, \\
\pi_2 & = \textstyle{\frac{1}{95040}}
       (17\th^4 + 1048\th^3 - 18486\th^2 - 271440\th + 1590381) \quad
       (\norm(\pi_2)=-2)
\end{align*}
\item
Prime factorization of $3$: 
\begin{align*}
3 & = -\pi_{31}\pi_{32} \\
\pi_{31} & = \textstyle{\frac{1}{95040}}
     (13\th^4 + 896\th^3 - 7806\th^2 - 280008\th - 861975) \quad
       (\norm(\pi_{31})=3)  \\
\pi_{32} & = \textstyle{\frac{1}{31680}} 
          (\th^4 + 68\th^3 - 810\th^2 - 22428\th + 105489) \quad
             (\norm(\pi_{32})=-3^4)   
\end{align*}
\item
Prime factorization of $5$: \label{page pi5}
\begin{align*}
5 & =  -\ep_1\ep_2^{-1}\ep_3^2\ep_4\pi_5^5 \\
\pi_5 & = \textstyle{\frac{1}{95040}}
             (13\th^4 + 896\th^3 - 7806\th^2 - 280008\th - 671895) \quad
             (\norm(\pi_5)=5)
\end{align*}
\item
Prime factorization of $11$: \label{page factor11}
\begin{align*}
11 & = \ep_1^{-1}\ep_2^{-1}\ep_4^{-1}\pi_{111}^2\pi_{112}^2\pi_{113} \\
\pi_{111} & = \textstyle{\frac{1}{47520}}
                (-\th^4 - 56\th^3 + 1554\th^2 + 22104\th + 43119) \quad
     (\norm(\pi_{111}=11) \\
\pi_{112} & = \textstyle{\frac{1}{95040}}
                  (\th^4 + 68\th^3 - 810\th^2 - 6588\th + 57969)  \quad
           (\norm(\pi_{112})=11) \\
\pi_{113} & = \textstyle{\frac{1}{31680}} 
              (-13\th^4 - 896\th^3 + 7806\th^2 + 280008\th + 798615)  \quad
              (\norm(\pi_{113})=-11)      
\end{align*}    \label{page def pi_113}
\end{itemize}
The above information combined with (\ref{eq ThueMahler5}) easily implies that
we have the following possibilities:
\begin{eqnarray}
\ideal{x-y\th} & = & \ideal{\pi_2}^5\ideal{\pi_{31}}^4\ideal{\pi_5}^{z_1}
       \ideal{\pi_{111}}^{w_1}\ideal{\pi_{112}}^{w_2}\ideal{\pi_{113}}^{w_3} 
     \label{eq alternative 1} \\
\ideal{x-y\th} & = &  \ideal{\pi_2}^5\ideal{\pi_{32}}\ideal{\pi_5}^{z_1}
       \ideal{\pi_{111}}^{w_1}\ideal{\pi_{112}}^{w_2}\ideal{\pi_{113}}^{w_3}  
       \label{eq alternative 2} 
\end{eqnarray} 
where, in both cases,
\begin{equation} \label{z2=w1+w2+w3}
w_1+w_2+w_3=z_2. 
\end{equation}
\subsubsection{Treating $5$-adically equations (\ref{eq alternative 1}) and (\ref{eq alternative 2}) }
                  \label{subsubsec working 5adically}
From  the ideal equations (\ref{eq alternative 1}) and (\ref{eq alternative 2}) we will compute an upper
bound for the unknown exponent $z_1$, using the ``Second Corollary of Lemma 1'' in Section 5 
of \cite{TdW}. As a consequence, the equation (\ref{eq ThueMahler5}) will be replaced by a rather small
number of similar equations in which only the exponent $z_2$ will be unknown; this is certainly a gain.
With the notations of Sections 3,5 of \cite{TdW} we have in our case:
%
%\pagebreak
\begin{center}
\tablecaption{Application of ``Second Corollary of Lemma 1'' in \cite{TdW} when $p=5$} 
 \label{Table Sections 3,5 TdW when p=5}  
\tablefirsthead{ 
\multicolumn{2}{c}{} \\[-8mm] \hline
   &  \\[-3mm]
Notations in \cite[\textsection\textsection\, 3,5]{TdW}  
                        & Interpretations in this case and [references in this paper]  \\ \hline\hline
               }
\tablehead{ \hline
\multicolumn{2}{|l|}{{\small \sl continued from previous page}} \\ \hline
Notations in \cite[\textsection\textsection\, 3,5]{TdW}  & Interpretations in this paper  \\
          }          
 \tabletail{\hline \multicolumn{2}{|r|}{\small\sl continued on next page}\\ \hline}
\tablelasttail{\hline\hline }
 \begin{supertabular}{|c|c|}
$p$ & $5$ \\ \hline
    &      \\[-4mm] 
$g(t)$\, (\textsection 3) & $t^5+65t^4-870t^3-18990t^2+26325t+255069\;\;$ [(\ref{eq def g})]  \\ \hline
    &      \\[-4mm] 
$m$ \, (\textsection 3) 
     & $1$; $g(t)$ is irreducible over $\Q_5$, hence $g(t)=g_1(t)$ \\ \hline
$\idp_1,\, e_1,\, d_1$\, (\textsection 3) 
               & $\ideal{\pi_5},\, 5,\, 1\quad $ [page \pageref{page pi5}] \\ \hline 
$e$ \, (\textsection 5) & $5$  \\ \hline               
    &      \\[-4mm] 
$D_{\th}$\, (\textsection 5) & $2^{32}3^{12}5^{11}11^6\quad $ [(\ref{eq def g})] \\ 
 \end{supertabular}
 \end{center}
In view of the fact that $m=1$ and $e_1=5$, the ``Second Corollary of Lemma 1'' in \cite{TdW}
implies that 
$ z_1=\ord{\pi_5}{x-y\th}\leq\frac{1}{2}e\cdot\ord{5}{D_{\th}}=55/2 $, hence,
\begin{equation} \label{eq bound z1}
z_1 \leq 27, \quad\mbox{implying $a\leq 55$}.
\end{equation}
\subsubsection{Treating $11$-adically equations (\ref{eq alternative 1}) and (\ref{eq alternative 2}) }
                  \label{subsubsec working 11adically}
Once again we use the notations of Sections 3,5 of \cite{TdW}. 
Now $g(t)=g_1(t)g_2(t)g_3(t)$ is the factorization of $g(t)$ into irreducible 
polynomials of $\Q_{11}[t]$, where
\[
g_1(t)= t^2+(3 + 3\cdot 11 + 8\cdot 11^2 + 9\cdot 11^3 + 5\cdot 11^4 + \cdots)t  
+(5 + 7\cdot 11^2 + 10\cdot 11^3 + 10\cdot 11^4 +\cdots),
\]
\[
g_2(t)=t^2+(3 + 5\cdot 11 + 5\cdot 11^2 + 2\cdot 11^4 + \cdots)t 
          + (5 + 11^3 + 7\cdot 11^4+\cdots),
\]          
\[
g_3(t) = t-(7 + 2\cdot 11 + 2\cdot 11^2 + 10\cdot 11^3 + 7\cdot 11^4 +\cdots).
\]
Let
$ g_1(\th_1)=0,\; g_2(\th_2)=0,\; g_3(\th_3)=0$.
Following the notation of the beginning of \textsection 5 of \cite{TdW} we 
denote the $\Q_{11}$-conjugates of the $\th_i$'s as follows:
\begin{itemize}
\item $\cnj{\th_1}{i}$, $i=1,2$; the roots of $g_1(t)$, living in a quadratic
extension of $\Q_{11}$.
\item $\cnj{\th_2}{i}$, $i=1,2$; the roots of $g_2(t)$, living in a quadratic
extension of $\Q_{11}$.
\item $\cnj{\th_3}{1} = \th_3=
7 + 2\cdot 11 + 2\cdot 11^2 + 10\cdot 11^3 + 7\cdot 11^4 +\cdots \in\Q_{11}$;
the root of $g_3(t)$.
\end{itemize}
%
%\pagebreak
\begin{center}
\tablecaption{Application of ``Second Corollary of Lemma 1'' in \cite{TdW}  when $p=11$} 
 \label{Table Sections 3,5 TdW when p=11}  
\tablefirsthead{ 
\multicolumn{2}{c}{} \\[-8mm] \hline
          &  \\[-3mm]
Notations in \cite[\textsection\textsection\, 3,5]{TdW}  
                        & Corresponding values in this case \\
                        &     and [references in this paper]
                                                             \\ \hline\hline
               }
\tablehead{ \hline
\multicolumn{2}{|l|}{{\small \sl continued from previous page}} \\ \hline
Notations in \cite[\textsection\textsection\, 3,5]{TdW}   & Interpretations and [references] in this paper  \\
          }                   
\tabletail{\hline \multicolumn{2}{|r|}{\small\sl continued on next page}\\ \hline}
\tablelasttail{\hline\hline}
 \begin{supertabular}{|c|c|}
$p$ & $11$ \\ \hline
   &      \\[-4mm] 
$g(t)$\, (\textsection 3) & $t^5+65t^4-870t^3-18990t^2+26325t+255069\;\;$  [(\ref{eq def g})]  \\ \hline
    &      \\[-4mm] 
$m$ \, (\textsection 3) 
     & $3$; $g(t)=g_1(t)g_2(t)g_3(t)\;\;$  [begining of Subsection \ref{subsubsec working 11adically}\\ \hline
       & \\[-10pt] 
 \hspace{-12mm}$\idp_1,\, e_1,\, d_1$
               &  $\ideal{\pi_{111}},\, 2,\, 1$ \\ 
$\idp_2,\, e_2,\, d_2$  \hspace*{3mm} (\textsection 3) 
               & $\ideal{\pi_{112}},\, 2,\, 1$ \\ 
      \hspace{-13mm} $\idp_3,\, e_3,\, d_3$  
               & $\ideal{\pi_{113}},\, 1,\, 1$ \\   
          & \hspace*{50mm} [page \pageref{page factor11}]  \\  \hline                          
$e$ \, (\textsection 5) & $1$  \\ \hline               
    &      \\[-4mm] 
$D_{\th}$\, (\textsection 5) & $2^{32}3^{12}5^{11}11^6\quad $ [(\ref{eq def g})]  \\
 \end{supertabular}
 \end{center}
Since we intend to apply the \emph{Prime Ideal Removing Lemma} 
\cite[Lemma 1]{TdW}, we must compute 
\begin{equation} \label{eq cond for PrIdRemLem}
\max\{e_i,e_j\}\cdot\ord{11}{\cnj{\th_i}{k}-\cnj{\th_j}{l}},
\quad i,j\in\{1,2,3\},\,i\neq j,
\end{equation}
where $k=1$ if $i=3$ and $k\in\{1,2\}$ if $i=1$ or $2$; and analogously,
$l=1$ if $j=3$ and $l\in\{1,2\}$ if $j=1$ or $2$.
According to the discussion in Appendix \ref{Append working padically}, in order to compute 
$\ord{11}{\cnj{\th_i}{k}-\cnj{\th_j}{l}}$ for fixed $1\leq i<j\leq 3$, it suffices to compute a polynomial 
$h_{ij}(t)\in\Q_{11}[t]$ such that $h_{ij}(\th_j-\th_i)=0$.
Moreover, since $\ord{p}{-\al}=\ord{p}{\al}$, it is clear that it suffices
to consider only the values $(i,j)=(1,3),(2,3),(1,2)$.
\\ \noindent
Obviously, $h_{i3}(t)=g_i(t+\th_3)$, hence
\[
h_{13}(t)= t^2+(6 + 8\cdot 11 + 11^2 + 8\cdot 11^3 + 10\cdot 11^4 + 
\cdots)t + (9 + 6\cdot 11 + 2\cdot 11^2 + 8\cdot 11^3 +\cdots).
\]
\[
h_{23}(t)=t^2+(6 + 10\cdot 11 + 9\cdot 11^2 + 9\cdot 11^3 + 6\cdot 11^4 
+ \cdots)t + (9 + 9\cdot 11 + 11^2 + 9\cdot 11^3 + 9\cdot 11^4 +\cdots),
\]
It follows that, when in (\ref{eq cond for PrIdRemLem}) we have
$(i,j)=(1,3),(2,3)$, then $\ord{11}{\cnj{\th_i}{k}-\cnj{\th_j}{l}}=0$.
\\ \noindent
When $(i,j)=(1,2)$, the following lemma, whose proof is a matter of straightforward calculations, 
gives us a quartic polynomial $h_{12}(t)\in\Q_{11}[t]$ which has $\cnj{\th_1}{k}-\cnj{\th_2}{l}$ 
as a zero (independent from $k,l$). 
\begin{lemma} \label{lemma th_i-th_j}
If $\th_i^2+a_i\th_i+b_i=0$ for $i=1,2$, then $\th_2-\th_1$ is a root of
$t^4+c_3t^3+c_2t^2+c_1t+c_0$, where
\[
c_3 =2(a_2-a_1),\quad c_2=a_2^2+a_1^2-3a_1a_2+2b_2+2b_1
\]
\[
c_1=a_1^2a_2-a_2^2a_1-2b_2a_1-2a_1b_1+2a_2b_2+2b_1a_2,\quad
\]
\[
c_0=b_2^2+b_1^2+b_2a_1^2+b_1a_2^2-b_2a_1a_2-b_1a_1a_2-2b_2b_1.
\]
\end{lemma}
The constant term of $h_{12}(t)$ is 
$9\cdot 11^2 + 5\cdot 11^3 + 6\cdot 11^4 + \cdots$, hence 
(\ref{eq order formula}) gives  
$\ord{11}{\cnj{\th_1}{k}-\cnj{\th_2}{l}}=2/4=1/2$.
\\ \noindent
In view of the above, when $(i,j)=(1,3)$ or $(2,3)$, the number 
(\ref{eq cond for PrIdRemLem}) is zero, therefore, statement (i) of the 
aforementioned Prime Ideal Removing Lemma  implies that $x-y\th$ is divisible 
by at most one prime among $\pi_{111}$ and $\pi_{113}$ (equivalently: 
$w_1=0$ or $w_3=0$) and by at most one prime among $\pi_{112}$ and $\pi_{113}$ 
(equivalently: $w_2=0$ or $w_3=0$), hence 
\begin{equation} \label{eq w3=0 or}
\mbox{either $w_3=0$ or $(w_1,w_2)=(0,0)$.}
\end{equation} 
When $(i,j)=(1,2)$, the number (\ref{eq cond for PrIdRemLem}) is equal to $1$,
hence, again by statement (i) of the Prime Ideal Removing Lemma, it follows that 
at most one among $\pi_{111}$ and $\pi_{112}$ divides $x-y\th$ with power $>1$. 
\emph{If this actually occurs} for $\pi_{11i}$ ($i=1$ or $2$), which means that
$\ord{\pi_{11i}}{x-y\th}>1$, then statement (ii) of the Prime Ideal Removing 
Lemma implies that
\[
\ord{\pi_{11i}}{x-y\th}\leq 2\,\ord{11}{\cnj{\th_i}{1}-\cnj{\th_i}{2}}
=2\cdot \frac{1}{2}=1,
\]
because, $(\cnj{\th_i}{1}-\cnj{\th_i}{2})^2$ being the discriminant of the
polynomial $g_i(t)$, is equal to either $6\cdot 11 + 8\cdot 11^2 + O(11^3)$
if $i=1$, or to $7\cdot 11 + 2\cdot 11^2 + O(11^3) $ if $i=2$.
This contradiction shows that $\ord{\pi_{11i}}{x-y\th}\leq 1$ for both 
$i=1,2$, i.e.
\begin{equation}  \label{eq w1,w2 at most 1}
w_i\leq 1 \quad (i=1,2).
\end{equation}
If we combine (\ref{eq w3=0 or}) and (\ref{eq w1,w2 at most 1}) we see
that we have the following possibilities:
\begin{equation} \label{eq possibilies for the w's}
(w_1,w_2,w_3)= (0,0,w_3),\,(0,1,0),\,(1,0,0),\,(1,1,0),
\end{equation}
where in the first case we understand that $w_3$ can be ``large''.
The remaining three possibilities combined with the relations
(\ref{eq alternative 1}) and (\ref{eq alternative 2}), lead us to
\begin{eqnarray*}
\ideal{x-y\th} & = & \ideal{\pi_2}^5\ideal{\pi_{31}}^4\ideal{\pi_5}^{z_1}
       \ideal{\pi_{111}}^{w_1}\ideal{\pi_{112}}^{w_2} 
     \\
\ideal{x-y\th} & = &  \ideal{\pi_2}^5\ideal{\pi_{32}}\ideal{\pi_5}^{z_1}
       \ideal{\pi_{111}}^{w_1}\ideal{\pi_{112}}^{w_2}.  
\end{eqnarray*} 
By (\ref{z2=w1+w2+w3}), $z_2=w_1+w_2+w_3=w_1+w_2=1,2$, and by (\ref{eq bound z1}),
$0\leq z_1\leq 27$. Taking norms in the above relations, we obtain the following fifty six Thue equations
(cf.~(\ref{eq ThueMahler5}):
\begin{eqnarray}
\lefteqn{x^5+65x^4y-870x^3y^2-18990x^2y^3+26325xy^4+255069y^5=c} \label{Thue eqs} \\
  & & c\in\{-2^5\cdot 3^4\cdot 5^{z_1}\cdot 11^{z_2}:\:  \nonumber
    0\leq z_1\leq 27,\, 1\leq z_2\leq 2 \}.
\end{eqnarray}
The {\sc magma} routine for solving Thue equations, based on Bilu \& Hanrot method \cite{BiHa}
(which improves the method of \cite{TdW-Thue}) ``answers'' that there are no solutions at all.
The computation cost for this task is less than 2.5 seconds.
\\ \noindent
In view of the above discussion, we are left with the first case 
in (\ref{eq possibilies for the w's}), hence we have to solve the ideal equations 
$ \ideal{x-y\th}  =  \ideal{\pi_2}^5\ideal{\pi_{31}}^4\ideal{\pi_5}^{z_1}
       \ideal{\pi_{113}}^{w_3}$
and       
$\ideal{x-y\th}  =  \ideal{\pi_2}^5\ideal{\pi_{32}}\ideal{\pi_5}^{z_1}
       \ideal{\pi_{113}}^{w_3}$, 
where, in both cases, $ 0\leq z_1\leq 27$.   

To sum up, the solution of the equation (\ref{eq ThueMahler5}) is reduced to that of the 
equation  
\begin{eqnarray} 
\lefteqn{x-y\th  =  
\al\ep_1^{a_1}\ep_2^{a_2}\ep_3^{a_3}\ep_4^{a_4}\pi_{113}^{n_1}
        }  \label{eq final Thue-Mahler} \\
    & & \al\in\{\pi_2^5\pi_{31}^4\pi_5^{z_1},\,\pi_2^5\pi_{32}\pi_5^{z_1}\,:
      \: 0\leq z_1\leq 27\},\quad n_1=w_3=z_2=(b-1)/2. \nonumber
\end{eqnarray}
in the unknowns $(a_1,a_2,a_3,a_4,n_1)\in\Z^4\times\Z_{\geq 0}$.
\subsubsection{From equation (\ref{eq final Thue-Mahler}) to $S$-unit equation (\ref{eq (12) of TdW}) }
                 \label{subsubsec S-unit eq}
Let $K$ be an extension of $F$ such that $g(t)$ has at least three linear factors in $K[t]$.
Actually, in our case, such an extension coincides with  the splitting field of $g(t)$ over $F$ 
(see (\ref{eq def g})). 
We have $K=\Q(\om)$ and the minimal polynomial of $\om$ over $\Q$, denoted by $G(t)$, is
of degree $20$ (see Appendix \ref{Append field K}). 
Thus, there exist $\cnj{\th}{i}(t)\in\Q[t]$  ($i=1,\ldots,5$), so that  the $\Q$-conjugates 
$\cnj{\th}{i}$ of $\th$ are
\[
\cnj{\th}{i}(\om)\in\Q(\om)=K   \qquad (i=1,\ldots,5).
\]
For every $i\in\{1,\ldots,5\}$, the $i$-th embedding $F\hookrightarrow K$ is
characterized by $\th\mapsto\cnj{\th}{i}(\om)$ and maps the general element 
$\be\in F$ to its $i$-th conjugate $\cnj{\be}{i}(\om)$. This belongs to $\Q(\om)$, hence it is a polynomial 
expression in $\om$, of degree at most $19$, with rational coefficients.

On the other hand, if $\idP$ is the prime ideal of $K$ over $\ideal{\pi_{113}}$, mentioned in 
Appendix \ref{Append field K},\footnote{See just above and below of relation (\ref{eq e,f of idP}).}
then, by the discussion of Appendix \ref{Append working padically}, there is an embedding  
$K\hookrightarrow K_{\idP}=\Q_{11}(\om_{\idP})$, where $G_{\idP}(\om_{\idP})=0$ for a specific 
second-degree factor $G_{\idP}(t)$ of $G(t)$, irreducible over $\Q_{11}$; 
see Appendix \ref{Append field K}.
This embedding is characterized by $\om\mapsto\om_{\idP}$, so that the $11$-adic roots of $g(t)$ are
\[
\cnj{\th}{i}(\om_{\idP})\in\Q(\om_{\idP})=K_{\idP} \qquad (i=1,\ldots,5)
\]
and, for every $\be\in F$, if the $i$-th conjugate of $\be$ over $\Q$ is $\cnj{\be}{i}(\om)$ (see a few lines
above), then the embedding $\om\mapsto\om_{\idP}$ maps $\be$ to $\cnj{\be}{i}(\om_{\idP})$.
\\ \noindent
If we work $p$-adically with $p=11$, then, by $\cnj{\th}{i},\cnj{\be}{i},\ldots $ we will understand 
$\cnj{\th}{i}(\om_{\idP})$, $\cnj{\be}{i}(\om_{\idP}),\ldots $; and if we work $p$-adically with 
$p=\mbox{\emph{infinite prime}}$,  by $\cnj{\th}{i},\cnj{\be}{i},\ldots $ we will understand 
$\cnj{\th}{i}(\om),\cnj{\be}{i}(\om),\ldots $. Our discussion below applies to both cases of $p$.

Applying the $i$-th embedding to the relation (\ref{eq final Thue-Mahler}) we obtain the $i$-th conjugate 
relation
\[
x-y\cnj{\th}{i}  =  
\cnj{\al}{i}{\cnj{\ep_1}{i}}^{a_1}{\cnj{\ep_2}{i}}^{a_2}{\cnj{\ep_3}{i}}^{a_3}{\cnj{\ep_4}{i}}^{a_4}
{\cnj{\pi_{113}}{i}}^{n_1}. 
\]
\label{page conjugates}
Then, for $i=i_0,j,k$, where  $i_0,j,k\in\{1,\ldots,5\}$ are any three distinct indices, we obtain three 
conjugate relations, analogous to the above. Eliminating $x,y$ from the these three relations we finally 
obtain (cf.~\cite[Section 7]{TdW})
\begin{equation} \label{eq (12) of TdW}
\la:=\de_1\left(\frac{\cnj{\pi_{113}}{k}}{\cnj{\pi_{113}}{j}}\right)^{\!\!n_1}
\prod_{i=1}^4\left(\frac{\cnj{\ep_i}{k}}{\cnj{\ep_i}{j}}\right)^{\!\!a_i}-1 =
\de_2\left(\frac{\cnj{\pi_{113}}{i_0}}{\cnj{\pi_{113}}{j}}\right)^{\!\!n_1}
\prod_{i=1}^4\left(\frac{\cnj{\ep_i}{i_0}}{\cnj{\ep_i}{j}}\right)^{\!\!a_i}, 
\end{equation}
where
\begin{equation} \label{eq deltas}
\de_1=\frac{\cnj{\th}{i_0}-\cnj{\th}{j}}{\cnj{\th}{i_0}-\cnj{\th}{k}}\cdot\frac{\cnj{\al}{k}}{\cnj{\al}{j}},
\quad
\de_2=\frac{\cnj{\th}{j}-\cnj{\th}{k}}{\cnj{\th}{k}-\cnj{\th}{i_0}}\cdot\frac{\cnj{\al}{i_0}}{\cnj{\al}{j}}.
\end{equation}
Now and until the end of the paper we put
\[
 H =\max\{n_1,\,|a_1|,\,|a_2|,\,|a_3|,\,|a_4|\}
\]
\subsubsection{Equation (\ref{eq (12) of TdW}) implies an upper bound $n_1\leq c_{13}\log H$ }
    \label{subsubsec n1<c13(logH+c14)}
We will prove the inequality in its title of this subsection, where $c_{13}$ is given by 
(\ref{N<c13(logH+c14)}). Our main tool is the important Theorem \ref{thm Yu thm simplified} due to 
Kunrui Yu which, given the algebraic numbers $\al_1,\ldots,\al_n$ and a prime $p$, provides 
an upper bound for the $p$-adic valuation of $\al_1^{b_1}\cdots\al_n^{b_n} - 1$, for any 
$b_1,\ldots,b_n\in\Z$,  in terms of $\log \max\{3, |b_1|,\ldots,|b_n|\}$.

We turn to the relation (\ref{eq (12) of TdW}), which we view as an algebraic relation over $\Q_{11}$. 
According to the discussion in Appendix \ref{Append field K}, the $11$-adic roots $\cnj{\th}{i}\in\C_{11}$  of 
$g(t)$ are identified with $\th_i(\om_{\idP})$ ($i=1,\ldots,5$). 

We choose the indices $i_0,j,k$ following the
instructions in \cite{TdW}, bottom of p.~235 and beginning of p.~236 up to Lemma 3. 
According to the discussion therein, since $\pi_{113}$ corresponds to the polynomial $g_3(t)$ whose root 
is $\th_5(\om_{\idP})$ (cf.~end of Appendix \ref{Append field K}), we must choose $i_0=5$; 
\label{page choice of i0,j,k}
and since $\th_1(\om_{\idP})$ and $\th_3(\om_{\idP})$ are (according to the end of 
Appendix \ref{Append field K}, again) roots of the quadratic irreducible polynomial $g_1(t)\in\Q_{11}[t]$, we 
can choose $j=1$ and $k=3$.
In view of \cite[Lemma 3\,(i)]{TdW}, $\ord{11}{  \cnj{\pi_{113}}{k}/\cnj{\pi_{113}}{j}}=0$ and
by  \cite[Corollary of Lemma 2\,(i)]{TdW},  $\ord{11}{  \cnj{\ep_i}{k}/\cnj{\ep_i}{j}}=0$ for $i=1,\ldots,4$.
Also, since $\cnj{\th}{k},\cnj{\th}{j}$ are $11$-adic roots of a second degree irreducible polynomial over 
$\Q_{11}$, it follows, according to the second ``bullet'' in page 236 of \cite{TdW}, that $\ord{11}{\de_1}=0$.
These facts will be used in the application of Theorem \ref{thm Yu thm simplified}.  

Also, the relation (13) of \cite[Theorem 5]{TdW} holds, which in our case reads
$\ord{11}{\la}=\ord{11}{\de_2}+n_1$ \footnote{Actually, according to the relation (13) of 
\cite[Theorem 5]{TdW}, $n_1$ is multiplied by a positive integer $h_1$, defined in 
\cite[Section 6]{TdW}, which is a divisor of the order of the ideal-class group. 
In our case, the ideal-class group is trivial, hence $h_1=1$.}.
A computation shows that $\ord{11}{\de_2}=1/2$
 \footnote{By (\ref{eq deltas}) and (\ref{eq final Thue-Mahler}) there are 56 possible values for $\de_2$.},
hence 
\begin{equation} \label{eq 1st ub for n1}
\ord{11}{\la}=n_1 +\textstyle{\frac{1}{2}}
\end{equation}
Now we are ready to apply Theorem \ref{thm Yu thm simplified}. With four minor corrections, this is 
Theorem 11.1 of K.~Hambrook's thesis \cite{Hamb}. It is a consequence of Theorems 1 and 3 of \cite{YuIII} 
and the Lemma in the Appendix of \cite{YuII}.
\begin{theorem}[Kunrui Yu]
           \label{thm Yu thm simplified}
 Let $\al_1,\ldots,\al_n$ ($n\geq 2$) be non-zero algebraic numbers and
\[
K=\Q(\al_1,\ldots,\al_n),\quad D=[K:\Q].
\]
Let $p$ be a rational prime, $\idP$ a prime ideal of the ring of integers of $K$ lying above $p$ and
$e_{\idP}=e_{K/\Q}(\idP)$, $f_{\idP}=f_{K/\Q}(\idP)$ the ramification index and residue class degree,
respectively, of $\idP$. 
\\
\noindent
Now define $d$ and $f$ as follows:
\\
If $p=2$ then
\[
d=\begin{cases}
     D & \mbox{if $e^{2\pi/3}\in K$} \\
     2D & \mbox{if $e^{2\pi/3}\not\in K$} 
    \end{cases},
    \qquad
 f=\begin{cases}
     f_{\idP} & \mbox{if $e^{2\pi/3}\in K$} \\
     \max\{2,f_{\idP}\} & \mbox{if $e^{2\pi/3}\not\in K$}
    \end{cases}.
\]
\hspace*{3mm}
If $p\geq 3$ and $p^{f_{\idP}}\equiv 3\pmod 4$ then 
\[
d=D,  \qquad f=f_{\idP}.
\]
If $p\geq 3$ and $p^{f_{\idP}}\equiv 1\pmod 4$ then 
\[
d=\begin{cases}
     D & \mbox{if $e^{2\pi/4}\in K$} \\
     2D & \mbox{if $e^{2\pi/4}\not\in K$} 
    \end{cases},
    \qquad
 f=\begin{cases}
     f_{\idP} & \mbox{if $e^{2\pi/4}\in K$ or $p\equiv 1\pmod 4$} \\
     \max\{2,f_{\idP}\} & \mbox{if $e^{2\pi/3}\not\in K$ and $p\equiv 3\pmod 4$}
    \end{cases}.
\]
Put
\[
\tau=\frac{p-1}{p-2},   \qquad 
\kappa=\left\lceil\log\left(\frac{2e_{\mathfrak{p}}}{p-1}\right)/\log p\right\rceil,
\qquad
Q=\begin{cases}
      3 & \mbox{if $p=2$} \\
      4 & \mbox{if $p\geq 3$ and $p^f\equiv 1\pmod 4$} \\
      1 & \mbox{if $p\geq 3$ and $p^f\equiv 3\pmod 4$}
    \end{cases},
\]
\[
(\kappa_1,\kappa_2,\kappa_3,\kappa_4,\kappa_5,\kappa_6)=
\begin{cases}
(160,32,40,276,16,8) & \mbox{if $p=2$} \\
(759,16,20,1074,8,4) & \mbox{if $p=3,\,d\geq 2$} \\
(537,16,20,532,8,4) & \mbox{if $p=3,\,d=1$} \\
(1473,8\tau,10,394\tau,8,4) & \mbox{if $p\geq 5,\,e_{\idP}=1,\,p\equiv 1\pmod 4$} \\
(1282,8\tau,10,366\tau,8,4) & \mbox{if $p\geq 5,\,e_{\idP}=1,\,p\equiv 3\pmod 4,\,d\geq 2$} \\
(1288,8\tau,10,396\tau,8,4) & \mbox{if $p\geq 5,\,e_{\idP}=1,\,p\equiv 3\pmod 4,\,d=1$} \\
(319,16,20,402,8,4) & \mbox{if $p=5,\,e_{\idP}\geq 2$} \\
(1502,16,20,1372,8,4) & \mbox{if $p\geq 7,\,e_{\idP}\geq 2,\,p\equiv 1\pmod 4$} \\
(2190,16,20,1890,8,4) & \mbox{if $p\geq 7,\,e_{\idP}\geq 2,\,p\equiv 3\pmod 4$}
\end{cases},
\]
\[
c_2=\frac{(n+1)^{n+2}d^{n+2}}{(n-1)!}\frac{p^f}{(f\log p)^3}\max\{1,\log d\}
\max\{\log(e^4(n+1)d),\,e_{\idP},\,f\log p\},
\]
\[
c'_3=\kappa_1\kappa_2^n\left(\frac{n}{f\log p}\right)^{n-1}
\prod_{i=1}^n \max\left\lbrace h(\al_i),\,\frac{f\log p}{\kappa_3(n+4)d}\right\rbrace,
\]
\[
c''_3=\kappa_4(e\kappa_5)^np^{(n-1)\kappa}
\prod_{i=1}^n \max\left\lbrace h(\al_i),\,\frac{1}{e^2\kappa_6p^{\kappa}d}\right\rbrace.
\]
Let $b_1,\ldots,b_n$ be rational integers and define
\[
\la=\al_1^{b_1}\cdots\al_n^{b_n} - 1, \quad B=\max\{3, |b_1|,\ldots,|b_n|\}.
\]
If $\la\neq 0$ and $\ord{\idP}{\al_i}=0$  for $i=1,\ldots,n$, then
\[
\ordp(\la) < c'_{10}\log B, 
\footnote{We use the notation $c'_{10}$ in order to conform with the notation of \cite[page 238]{TdW}.}
\quad  c'_{10}=\frac{c_2\min\{c'_3,c''_3\}}{Q\cdot e_{\idP}}.
\]
\end{theorem}
Now we apply Theorem \ref{thm Yu thm simplified} to the $\la$ given in (\ref{eq (12) of TdW}), as
interpreted in the beginning of this section, with $i_0=5,j=1, k=3$. Our application is briefly described in
\mbox{Table \ref{Table Yu}.}
\\
%
%\begin{minipage}{\textwidth}
\begin{center}
\tablecaption{Application of Theorem \ref{thm Yu thm simplified}} 
 \label{Table Yu}  
\tablefirsthead{ 
\multicolumn{2}{c}{} \\[-7mm] \hline
  &      \\[-10pt]
Notations in Theorem \ref{thm Yu thm simplified}    & Values in this paper
                                                             \\ \hline\hline
               }
\tablehead{ \hline
\multicolumn{2}{|l|}{{\small \sl continued from previous page}} \\ \hline
Notations in Theorem \ref{thm Yu thm simplified}   & Values  in this paper  \\
          }          
 \tabletail{\hline \multicolumn{2}{|r|}{\small\sl continued on next page}\\ \hline}
\tablelasttail{\hline\hline }          
\tabletail{\hline \multicolumn{2}{|r|}{\small\sl continued on next page}\\ \hline}
\tablelasttail{\hline\hline}
 \begin{supertabular}{|c|c|}
      &           \\[-10pt]
$n$ & $6$ \\ \hline
$(\al_1,b_1)$ & $(\de_1,1)$   [equation (\ref{eq deltas}) with $(i_0,j,k)=(5,1,3)$]   \\ \hline
        & \\[-10pt]
  $(\al_2,b_2)$   & $(\cnj{\pi_{113}}{3}/\cnj{\pi_{113}}{1},n_1)$ \;  [equation (\ref{eq (12) of TdW})] \\ \hline    
       & \\[-10pt]
$(\al_i,b_i)$, ($i=3,4,5,6$) 
     & $(\cnj{\ep_i}{3}/\cnj{\ep_i}{1},a_{i-2})$, ($i=3,4,5,6$)\;
          [equation (\ref{eq (12) of TdW})] \\ \hline
            & \\[-10pt]
$K$ &  $K$\; [Appendix \ref{Append field K}] \\ \hline
    &      \\[-4mm] 
$D$ & $20$ \\ \hline
$p$  & $11$ \\ \hline
       & \\[-10pt]
$\idP$ & $\idP$\; [Appendix \ref{Append field K} just above equation (\ref{eq e,f of idP})] \\ \hline
        & \\[-10pt]
$(f_{\idP},\,e_{\idP})=(f_{K/\Q}(\idP),\,e_{K/\Q}(\idP))$ & $(1,2)$; [(\ref{eq e,f of idP})] \\ \hline
$B$ & $\max\{3,n_1,|a_1|,|a_2|,|a_3|,|a_4|\}$ \\ \hline
           & \\[-10pt]
$c'_{10}$  & $< 9.9\cdot 10^{30}$ \\
 \end{supertabular}
  \end{center}  
%  \end{minipage}
%  \\[4mm]
 %
 A remark has its place here: 
By (\ref{eq final Thue-Mahler}) and the definition of $\de_1$ in (\ref{eq deltas}) we see that $\de_1$ runs 
through a set of cardinality 28, therefore, for each value of $\de_1$, we must compute the parameter 
$c'_{10}$. It turns out that, in all cases, $c'_{10}  < 9.9\cdot 10^{30}$ and this is mentioned in the above 
table.
Also, in the notation of \cite{TdW}, the use of Theorem \ref{thm Yu thm simplified} always implies 
$c'_{11}=0$. 

By writing a number of rather simple routines we automated the computations. Finally, by setting  
$c'_{10}\leftarrow \max_{\de_1} c'_{10}$ and  $c'_{11}\leftarrow \max_{\de_1} c'_{11}$ we find
$c'_{10}=9.99\cdot 10^{30}$ and $c'_{11}=0$.

%%%%%%%%%%%%%% mexri edw

By \cite[relation (14)]{TdW}, $n_1\leq c_{13}(\log H+c_{14})$, where $c_{13},c_{14}$ are explicitly 
computed from $c'_{10}$ and $c'_{11}$ following the simple instructions found on p.~238 of \cite{TdW}. 
The difference between the pairs $(c_{13},c_{14})$ and $(c'_{10},c'_{11})$, if any at all, is negligible in
practice. Anyway, in our case, it turns out easily that the two pairs coincide and, therefore,
\begin{equation} \label{N<c13(logH+c14)}
n_1 \leq c_{13} (\log H +c_{14}), \quad c_{13}=9.99\cdot 10^{30}, \; c_{14}=0,
\end{equation}
where
\begin{equation}  \label{def H}
  H =\max\{n_1,\,|a_1|,\,|a_2|,\,|a_3|,\,|a_4|\}.
\end{equation}
\emph{A computational remark}. According to the instruction of \cite[page 238]{TdW}, in order to compute 
$c_{13}$ from $c'_{10}$ we need the least positive integer $h$ such that $\idp^h$ is principal. In our case 
$\idp$ is already a principal ideal, therefore we take $h=1$. 
In order to compute $c_{14}$ from $c'_{11}$ we need to compute 
$\ord{11}{\de_2}$ for the $56$ values $\de_2$ (cf.~(\ref{eq (12) of TdW}) and (\ref{eq deltas})).
One shouldn't expect difficulties in carrying out such computations using {\sc magma} or any
other package specialized to Number Theory. 
\subsubsection{First explicit bounds for $H=\max\{n_1,|a_1|,|a_2|,|a_3|,|a_4|\}$ and $n_1$}
                                                               \label{subsubsec 1st ub}
We will prove the numerical upper bound (\ref{1st bound H}) for $H$,  based to  E.M.~Matveev's  lower 
bound for linear form in (real/complex) logarithms of algebraic numbers; see Theorem \ref{thm Matv} below.   
Then, as a straightforward consequence of (\ref{N<c13(logH+c14)}), this will imply the numerical
upper bound (\ref{1st bound n1}) for $n_1$.

We focus our attention to $1+\la$, where $\la$ is defined in relation (\ref{eq (12) of TdW}).
In this section we view $K$ embedded in the complex field $\C$, so that the algebraic
numbers appearing in $\la$ are complex numbers; actually, they are all \emph{real} numbers, because
all roots of $g(t)$ are real.  
Note that the indices $i_0,j,k$ figuring in (\ref{eq (12) of TdW}) are \emph{any} distinct indices from the set 
$\{1,\ldots,5\}$.
We follow step by step the very explicit instructions of Sections 9 and 10 of \cite{TdW} in order to compute
a chain of constants (in the order that are displayed below)
\[
c_{15},\:  c_{16}=0.129,\footnote{We give the value of $c_{16}$, because this will play a role later.}
\: c'_{17},\: c'_{18},\: c''_{17},\: c''_{18},\: c_{17},\: c_{18},\: c_{19},\: c_{20},\: 
c_{12},\: c_{21},\: c_{22}.
\]
This is a rather boring and cumbersome task if one performs the computations ``by hand'' (with the aid of a 
pocket calculator). Fortunately, the instructions are programmable in {\sc magma} without much difficulty,
so that the chain of computations is performed automatically. It turns out that $c_{22}=14$. 
According to the terminology of page 243 of \cite{TdW}, we are treating a ``real case''. Moreover, by
page 244 of \cite{TdW}, if we assume that $H>c_{22}=14$ ($H$ is defined in (\ref{def H})), 
then $1+\la$ is a positive real number and 
\begin{equation}   \label{eq log(1+la)}
\La=\log(1+\la) = 
\log|\de_1|+n_1\log\left|\frac{\cnj{\pi_{113}}{k}}{\cnj{\pi_{113}}{j}}\right|
+\sum_{i=1}^4a_i\log\left|\frac{\cnj{\al_{113}}{k}}{\cnj{\al_{113}}{j}}\right|.
\end{equation}
By a strong and handy result of E.M.~Matveev we can compute explicit constants
$c_7,c_8$ such that $\log (1+\la)>\exp (-c_7(\log H+c_8))$. 
More specifically we have the Theorem \ref{thm Matv} below, which is a slight restatement of 
Theorem 2.1 of \cite{Matv}. 
In this theorem $\log$ denotes an arbitrary but fixed branch of the logarithmic function on $\C$; if $x$ is a 
positive real number, $\log x$ always means real (natural) logarithm of $x$.  
\begin{theorem}\emph{(\cite[Theorem 2.1]{Matv})}
             \label{thm Matv}
Let $\La=b_1\log\al_1+\cdots +b_n\log\al_n$, where $b_1,\ldots,b_n\in\Z$ with $b_n\neq 0$, 
and $\al_1,\ldots,\al_n$ are algebraic numbers of degree at most $D$, embedded in $\C$, and 
$\log\al_1,\ldots,\log\al_n$ are linearly independent over $\Z$.
\\ \noindent
Consider $A_1,\ldots,A_n$ satisfying 
\[
A_i \geq \max\{Dh(\al_i),\: |\log\al_i|\} \qquad 1\leq i\leq n,
\]
where, in general, $h(\al)$ denotes the absolute logarithmic height of the algebraic number $\al$.
Set $\kappa=1$ if all $\al_i$'s are real; otherwise set $\kappa=2$. Next, define
\[ A=\max_{1\leq i\leq n}A_i/A_n,\quad   \Omega=A_1\cdots A_n \]
and
\[
B =\max_{1\leq i\leq n} |b_i|.
\]
Then 
\[
|\La|>\exp(-c_7(\log B +c_8)), 
\]
where
\[
c_7 =  \frac{16}{n!\kappa}e^n(2n+1+2\kappa)(n+2)(4(n+1))^{n+1}(en/2)^{\kappa} 
             \log(e^{4.4n+7}n^{5.5}D^2\log(eD))D^2\Omega,
\]
\[             
c_8 = \log(1.5eD\log(eD)A).
\]
\end{theorem}
\hspace*{\parindent}
Now we apply Theorem \ref{thm Matv} to the linear form $\La=\log(1+\la)$ in (\ref{eq log(1+la)}).
Following the instructions of \cite{TdW} (bottom of page 249 - beginning of page 250), we must consider 
$\La$ for all $i_0\in\{1,\ldots,5\}$, and for each specific $i_0$, the choice of the indices $j,k$ is
arbitrary, provided that $i_0\neq j\neq k\neq i_0$. 
 \label{page i0,j,k}
 Note that the condition of $\Z$-linear independence of the $\al_i$'s, imposed by Theorem \ref{thm Matv},  
in our case reads $\log(1+\la)\neq 0$. This is equivalent to $\la\neq 0$; we see that this is true by viewing
$\la$ as the right-hand side of the relation (\ref{eq (12) of TdW}).
The application of Theorem \ref{thm Matv} in our case is briefly described in Table \ref{Table Matv}.
In this table, $(i_0,j,k)$ runs through the set $\{(1,2,3),\,(2,1,3),\,(3,1,2),\,(4,1,2),\,(5,1,2)\}$.
We note that, the condition of $\Z$-linear independence of the $\al_i$'s, imposed by 
Theorem \ref{thm Matv}, in our case reads $\log(1+\la)\neq 0$. This is equivalent to $\la\neq 0$; 
we see that this is true by viewing $\la$ as the right-hand side of the relation (\ref{eq (12) of TdW}).
\begin{minipage}{\textwidth}
\begin{center}
\tablecaption{Application of Theorem \ref{thm Matv} to $\log(1+\la)$}
                 \label{Table Matv}  
\tablefirsthead{ 
\multicolumn{2}{c}{} \\[-7mm] \hline
 &      \\[-10pt]
 Notations in Theorem \ref{thm Matv}   
                        & Values in this paper    
                                                             \\ \hline\hline
               }
\tablehead{ \hline
\multicolumn{2}{|l|}{{\small \sl continued from previous page}} \\ \hline
Notations in Theorem \ref{thm Matv}   
                        & Values  in this paper
                                                             \\
          }
\tabletail{\hline \multicolumn{2}{|r|}{\small\sl continued on next page}\\ \hline}
\tablelasttail{\hline\hline}
 \begin{supertabular}{|c|c|}
      &           \\[-10pt]
$n$ & $6$ \\ \hline
$(\al_1,b_1)$ & $(\de_1,1)$   [equation (\ref{eq deltas}) with $(i_0,j,k)=(5,1,3)$]   \\ \hline
        & \\[-10pt]
  $(\al_2,b_2)$   & $(\cnj{\pi_{113}}{3}/\cnj{\pi_{113}}{1},n_1)$ \;  [equation (\ref{eq (12) of TdW})] \\ \hline    
       & \\[-10pt]
$(\al_i,b_i)$, ($i=3,4,5,6$) 
     & $(\cnj{\ep_i}{3}/\cnj{\ep_i}{1},a_{i-2})$, ($i=3,4,5,6$)\;
          [equation (\ref{eq (12) of TdW})] \\ \hline
            & \\[-10pt]
$K$ &  $K$\; [Appendix \ref{Append field K}] \\ \hline
    &      \\[-4mm] 
$D$ & $20$ \\ \hline
$p$  & $11$ \\ \hline
  $\idP$ & $\idP$\; [Appendix \ref{Append field K} just above equation (\ref{eq e,f of idP})] \\ \hline
        & \\[-10pt]
$(f_{\idP},\,e_{\idP})=(f_{K/\Q}(\idP),\,e_{K/\Q}(\idP))$ & $(1,2)$; [(\ref{eq e,f of idP})] \\ \hline
$B$ & $\max\{3,n_1,|a_1|,|a_2|,|a_3|,|a_4|\}$ \\ \hline
           & \\[-10pt]
$c_7,\; c_8$  & $<4.8626\cdot 10^{27},\;\; <5.7864$  \\ 
 \end{supertabular}
 \end{center}
 \end{minipage}
 \\[4mm]
We remark at this point that, actually, the values of $c_7,c_8$ which we obtain for the various choices
 of $(i_0,j,k)$ differ ``very little'', if they differ at all.

We\footnote{``We'' means ``our {\sc magma} code''.} continue to follow the instructions from the 
relation (24) of \cite{TdW} onwards and compute constants $c_{23},c_{24}, c_{25}$ and, finally, constants 
$c_{\mbox{{\footnotesize real}}}$ and $c_{27}$, such that  
\\
\textbullet\; $H=\max\{n_1,|a_1|,|a_2|,|a_3|,|a_4|\} <c_{\mbox{{\footnotesize real}}}\quad $ 
(see \cite[Theorem 10]{TdW}), 
\\
\textbullet\; $|\La|\leq c_{27}\exp\{-c_{16}\max_{1\leq i\leq 4}|a_i| \}\quad$
(see \cite[relation (29)]{TdW}).
\\ \noindent
Also, by \cite[Corollary of Theorem 10]{TdW}, 
$n_1\leq c_{13}(\log c_{\mbox{{\footnotesize real}}}+c_{14})$ and, in general, $n_1$ is considerably 
smaller than $c_{\mbox{{\footnotesize real}}}$.
\\ \noindent
According to our computations, the maximum value for $c_{\mbox{{\footnotesize real}}}$ \footnote{As 
$(i_0,j,k)$ runs through the set $\{(1,2,3),\,(2,1,3),\,(3,1,2),\,(4,1,2),\,(5,1,2)\}$.} is 
$1.3216\cdots\times 10^{43}$, from which we conclude that
\begin{equation} \label{1st bound H}
H < 1.3217\cdot 10^{43} =: K_0
\end{equation}  
and
\begin{equation} \label{1st bound n1}
n_1 <  9.918312\cdot 10^{32} =: N_0.
\end{equation} 
Also, $c_{27}<3.906653$; this constant, along with $c_{16}=0.129$ will be used in 
Subsection \ref{subsubsec padic reduc}.

{\bf Computation time}. The totality of computations that led to the bounds (\ref{1st bound H}) and
(\ref{1st bound n1}) was about 8 minutes.
\subsubsection{The first $p$-adic reduction}
                                 \label{subsubsec padic reduc}
In this section we reduce the upper bound (\ref{1st bound n1}) by a process we call \emph{$p$-adic 
reduction}, with $p=11$ in our case. 
imeFor the basic facts we refer to \cite{dW} and \cite[Sections 12,14,15]{TdW}. 
Given a rational prime $p$ and a $p$-adic number $x$ (in general, $x$ belongs to a finite extension of 
$\Q_p$), the $p$-adic logarithm of $x$ is denoted by $\log_p x$ and belongs to the same extension of 
$\Q_p$ in which $x$ belongs.

We go back to the relations (\ref{eq (12) of TdW}) and (\ref{eq deltas}). According to the discussion
in Appendix \ref{Append field K} and, more specifically, the notations etc on page 
\pageref{page choice of mathfrakP}, 
we have an embedding $K\hookrightarrow K_{\idP}$, where $K_{\idP}=\Q_{11}(\om_{\idP})$ is a quadratic 
extension of $\Q_{11}$ defined by the polynomial 
$G_{\idP}(t)=t^2+(10744341441+O(11^{10}))\,t+(9625552201+O(11^{10}))=0$,\footnote{see 
Table \ref{Table factorization in K} et.~seq. } which allows us to view the $\cnj{\th}{i}$'s figuring in the 
above relations as elements of $K_{\idP}$. According to our choice for $i_0,j,k$, made  in 
page \pageref{page choice of i0,j,k}, $(i_0,j,k)=(5,1,3)$ and, consequently,
$\cnj{\th}{i_0}=\cnj{\th}{5}=7050162550+O(11^{10})$, 
$\cnj{\th}{j}=\cnj{\th}{1}=(9038034724+O(11^{10}))\om_{\idP}+(8245826831+O(11^{10}))$,
$\cnj{\th}{k}=\cnj{\th}{3}=(4757114675+O(11^{10}))\om_{\idP}+(5113588460+O(11^{10}))$.
In the notation of page \pageref{page choice of i0,j,k}, these are the $11$-adic roots $\th_5(\om_{\idP})$,
$\th_1(\om_{\idP})$ and $\th_3(\om_{\idP})$ of $g(t)$, respectively\footnote{The two remaining $11$-adic 
 roots of $g(t)$ are 
$\th_2(\om_{\idP})=(517324682+O(11^{10}))\om_{\idP}+(5431351847+O(11^{10}))$ and
$\th_4(\om_{\idP})=(2621442663+O(11^{10}))\om_{\idP}+(7443205770+O(11^{10}))$.}.
By (\ref{eq 1st ub for n1}),  $\ord{11}{\la}=n_1+\frac{1}{2}$, therefore, by \cite[Lemma 12]{TdW},
$\ord{11}{\La_{11}}=n_1+\frac{1}{2}$, where
\[
\La_{11}=\log_{11}\de_1+n_1\log_{11}\left(\frac{\cnj{\pi_{113}}{3}}{\cnj{\pi_{113}}{1}}\right)+
\sum_{i=1}^4a_i\log_{11}\left(\frac{\cnj{\ep_i}{3}}{\cnj{\ep_i}{1}}\right).
\]
Because $\de_1$ depends also on the choice of $\al$ (cf.~relations (\ref{eq deltas}) and 
(\ref{eq final Thue-Mahler})), there are 56 possibilities for $\La_{11}$. Therefore, in what follows, we 
assume that, having chosen 
$\al\in\{\pi_2^5\pi_{31}^4\pi_5^{z_1},\,\pi_2^5\pi_{32}\pi_5^{z_1}\,: \: 0\leq z_1\leq 27\}$,
we compute the 11-adic logarithms appearing in $\La_{11}$; except for $\log_{11}\de_1$,
the remaining logarithms are independent from $\al$. 

Note that the values of $\log_{11}$ above belong to $K_{\idP}$ and, therefore, they are of the 
form $x_0+x_1\om_{\idP}$, where $x_0,x_1\in\Q_p$. 
If we put 
\[
\log_{11}\de_1=\rho_0+\rho_1\om_{\idP},\; 
\log_{11}(\cnj{\pi_{113}}{3}/\cnj{\pi_{113}}{1})=\la_0+\la_1\om_{\idP},\;
\log_{11}(\cnj{\ep_i}{3}/\cnj{\ep_i}{1})=\mu_{i0}+\mu_{i1}\om_{\idP} 
\]
($i=1,2,3,4$), then $\La_{11}=\La_{11,0}+\La_{11,1}\om_{\idP}$, where
\[
\La_{11,0}=\rho_0+n_1\la_0+\sum_{i=1}^{4}a_i\mu_{i0}, \quad
\La_{11,1}=\rho_1+n_1\la_1+\sum_{i=1}^{4}a_i\mu_{i1}.
\]  
Following the instructions of \cite[p.p.~256-257]{TdW} we put for $i=0,1$:
\[
\La_{11,i}'=-\be_{0i}-n_1\be_{1i}-a_1\be_{2i}-a_3\be_{3i}+a_4,
\]
where,
\[
\be_{0,i}=-\rho_i/\mu_{4i},\quad   \be_{1i}=-\la_i/\mu_{4i}, \quad
 \be_{ji}=-\mu_{j-1,i}/\mu_{4i} \;\; (j=2,3,4).
\]
We divided by $\ord{11}{\mu_{4i}}$, because 
$\ord{11}{\mu_{4i}}\leq\min\{\ord{11}{\la_1},\ord{11}{\mu_{1i}},\ldots,\ord{11}{\mu_{4i}}\}$.

Following the detailed instructions of \cite[Section 15]{TdW}, we put $W=\lceil K_0/N_0\rceil$,
and we choose appropriately a number $\kappa >1$ --this will become clear below-- and an integer $m$ 
such that 
\begin{equation}   \label{eq choose m}
11^mW=\kappa K_0^5, \footnote{The exponent $5$ is equal to the 
number of the unknown exponents in (\ref{eq (12) of TdW}).} 
\end{equation}
Then, for $i=0,1$, we consider the lattice $\Ga_{mi}$  which is generated by the column-vectors of the 
matrix
\[
\left(\begin{array}{ccccc}
W & 0 & 0 & 0 & 0 \\
0 & 1 & 0 & 0 & 0 \\
0 & 0 & 1 & 0 & 0 \\
0 & 0 & 0 & 1 & 0 \\
\cnj{\be_{1i}}{m} & \cnj{\be_{2i}}{m} & \cnj{\be_{3i}}{m} & \cnj{\be_{4i}}{m} & 11^m
\end{array}\right),
\]
where, in general, for $\be\in\Q_{11}$, we denote by $\cnj{\be}{m}$ the integer of the interval 
$[0,\,11^m-1]$ for which $\ord{11}{\be-\cnj{\be}{m}}\geq m$.
We also consider the column vector
\[
\mathbf{y}_i=\left(\begin{array}{c} 0 \\ 0 \\ 0 \\ 0 \\ -\cnj{\be_{0i}}{m} \end{array}\right).
\]
Note that, in view of our remark after the definition of $\La_{11}$, there are 56 possible values for
the vector $\mathbf{y}_i$, but the lattices $\Ga_{mi}$ are independent from $\al$.
\\
\noindent
Let $\mathbf{c}_1,\ldots, \mathbf{c}_5$ be the column-vectors of an (ordered) LLL-reduced basis of 
$\Ga_{mi}$ and $s_1,\ldots,s_5\in\Q$ be such that $\mathbf{y}_i=\sum_{j=1}^5 s_j\mathbf{c}_5$.       
 Let $j_0$ be the maximum index $j\in\{1,\ldots,5\}$ for which $s_j\not\in\Z$ and denote by $\ndist{s_{j0}}$ 
 the distance of $s_{j0}$ from the nearest to it integer. Finally, put
\[
\ell(\Ga_{mi},\mathbf{y}_i) \geq
\begin{cases}
    \textstyle{\frac{1}{4}}|\mathbf{c}_1| & \mbox{if $\mathbf{y}_i=\mathbf{0}$} \\
     \textstyle{\frac{1}{4}}\ndist{s_{j0}}\!\cdot|\mathbf{c}_1| & \mbox{if $\mathbf{y}_i\neq\mathbf{0}$} 
  \end{cases}   
\]
No we apply \cite[Proposition 15]{TdW}, which, in our case reads: 
\begin{quote}
\emph{If 
\begin{equation}  \label{eq cond of padic reduc}
\ell(\Ga_{mi},\mathbf{y}_i) >\sqrt{WN_0^2+4K_0^2},
\end{equation}
then $n_1<m$.}\footnote{Actually, in accordance to \cite[Proposition 15]{TdW}, the upper bound for $n_1$ 
is $(m-l)/h$. By the fact that, over $11$, the prime ideals of $F$ are principal, and the definition of $h$ in 
\cite[page 234]{TdW}, we have $h=1$. Also, $l$ is a small number, explicitly determined in page 257 of 
\cite{TdW} and, more specifically, below the relation (32); in our case, it turns out that $l=0$.}
\end{quote}
Heuristically, when $\kappa$ --and, accordingly by (\ref{eq choose m}), also $m$-- are sufficiently large, 
it is ``reasonable'' to expect that condition (\ref{eq cond of padic reduc}) is satisfied, which would imply
an upper bound for $n_1$.
Choosing in (\ref{eq choose m}) $\kappa=100$, so that $m=206$, we check that the condition 
(\ref{eq cond of padic reduc}) is satisfied, for either $i=0$ or $i=1$, for all but 10 values of  $\mathbf{y}_i$;
for the ten exceptional values of $\mathbf{y}_i$ we take $\kappa=1000$, so that $m=207$, and then
(\ref{eq cond of padic reduc}) is satisfied. 
\\[3pt]
\emph{As a consequence, we conclude that $n_1\leq 207$.}
\\
\noindent
{\bf Computation time}. The computation cost for this reduction step was less than 1 minute.
\subsubsection{The first reduction over $\R$}
                                 \label{subsubsec real reduc}
We have the upper bound $K_0=1.32171\cdot 10^{43}$ for $H=\max\{n_1,|a_1|,|a_2|,|a_3|,|a_4|\}$
and, by the conclusion of Subsection \ref{subsubsec padic reduc}, we already know that 
$n_1\leq 207=:N_1$.
Thus, in (\ref{eq log(1+la)}), coefficients $n_1,a_1,\ldots,a_4$ of the linear form  $\log(1+\la)$  satisfy 
$n_1\leq 207$ and $\max_i |a_i|\leq K_0$. Referring to (\ref{eq log(1+la)}),
let as put 
\[
\La=\log(1+\la)=\rho +n_1\la_1+\sum_{i=1}^4 a_i\mu_i,
\]
where the meaning of the real numbers $\rho,\la_1,\mu_1\ldots,\mu_4$ is obvious. Once again we stress
the fact that these six real numbers depend on the choice of the indices $(i_0,j,k)$ 
(cf.~page \pageref{page i0,j,k}), and $\rho=\log|\de_1|$ (cf.~(\ref{eq log(1+la)}))
depends also on the choice of $\al$ (cf.~relations (\ref{eq deltas}) and (\ref{eq final Thue-Mahler})).
Since there are 5 choices for $(i_0,j,k)$ and 56 choices for $\al$, this implies that there are
$5\times 56=280$ possibilities for the linear form $\La$.  Therefore, in what follows, we assume that,
having chosen $\al\in\{\pi_2^5\pi_{31}^4\pi_5^{z_1},\,\pi_2^5\pi_{32}\pi_5^{z_1}\,: \: 0\leq z_1\leq 27\}$
and $(i_0,j,k)\in\{(1,2,3),\,(2,1,3),\,(3,1,2),\,(4,1,2),\,(5,1,2)\}$, we compute the real numbers
$\rho,\la_1,\mu_1\ldots,\mu_4$.

We follow the reduction process of \cite{dW}, as presented in \cite[Section 16]{TdW}. We put  
$W'=\lceil K_0/N_1\rceil$ --this is independent from the above choices-- and choose a
number $\kappa>1$ and an integer $C$ so that $CW'\approx\kappa K_0^5$.\footnote{As in Subsection 
\ref{subsubsec padic reduc}, the exponent $5$ is the number of the unknowns exponents in 
(\ref{eq (12) of TdW}).} How we choose $\kappa$ will become clear below; as it turns out in practice, 
$\kappa$ depends on $\al$ and $(i_0,j,k)$.
We consider the lattice $\Ga_C$  which is generated by the column-vectors of the matrix
\[
\left(\begin{array}{ccccc}
W' & 0 & 0 & 0 & 0 \\
0 & 1 & 0 & 0 & 0 \\
0 & 0 & 1 & 0 & 0 \\
0 & 0 & 0 & 1 & 0 \\
\phi_1 & \psi_1 & \psi_2 & \psi_3 & \psi_4
\end{array}\right),
\]
where $\phi_1=\lfloor C\la_1\rfloor, \psi_i=\lfloor C\mu_i\rfloor$ ($i=1,\ldots,4$). Also, we put 
$\phi_0=\lfloor C\rho\rfloor$ and consider the column-vector
\[
\mathbf{y}=\left(\begin{array}{c} 0 \\ 0 \\ 0 \\ 0 \\ -\phi_0 \end{array}\right).
\]
As in the previous section, we compute an (ordered) LLL-reduced basis of $\Ga_C$, say 
$\mathbf{c}_1,\ldots, \mathbf{c}_5$. Let $s_1,\ldots,s_5\in\Q$ be 
the coefficients of $\mathbf{y}$ with respect to this basis, denote by $j_0$ the maximum index 
$j\in\{1,\ldots,5\}$ for which $s_j\not\in\Z$ and by $\ndist{s_{j0}}$ 
 the distance of $s_{j0}$ from the nearest to it integer. Finally, put
\[
\ell(\Ga_C,\mathbf{y}) \geq
\begin{cases}
    \textstyle{\frac{1}{4}}|\mathbf{c}_1| & \mbox{if $\mathbf{y}=\mathbf{0}$} \\
     \textstyle{\frac{1}{4}}\ndist{s_{j0}}\!\cdot|\mathbf{c}_1| & \mbox{if $\mathbf{y}\neq\mathbf{0}$}. 
  \end{cases}   
\]
Following the instructions of \cite[page 265]{TdW} we put $R=N_1+4K_0+1$ and 
$S={W'}^2N_1^2+3K_0^2$. By \cite[Proposition 16]{TdW}: 
\begin{quote}
\emph{If $\;\ell(\Ga_C,\mathbf{y}) \geq \sqrt{R^2+S}\quad (*) $ 
\[
\mbox{then} \quad H \leq\frac{1}{c_{16}}\{\log c_{27}+\log C-\log(\sqrt{\ell(\Ga_C,\mathbf{y})^2-S}-R)\}.
\quad (**)
\]
}
\end{quote}
Heuristically, one can argue that, if $\kappa$ is sufficiently large and $C=\lceil\kappa K_0^5/W'\rceil$,
then it is ``reasonable'' to expect that $\ell(\Ga_C,\mathbf{y}) \geq \sqrt{R^2+S}$ and, consequently, an
upper bound for $H$ is obtained from $(**)$, which is of the size of $\log K_0$.
\\ \noindent
To give an idea, if $(i_0,j,k)=(1,2,3)$ and we take $\kappa=100$, $C=10^{187}$, then the condition
$(*)$ is satisfied for all $\al$'s and, as $\al$ runs through all its possible values, the maximum bound
$(**)$ is 229. 
If $(i_0,j,k)=(2,1,3)$ and $\kappa=100$, $C=10^{187}$, then $(*)$ holds for all but 11 values of $\al$. For 
the 45 ``successful'' values of $\al$ the maximum bound $(**)$ is 229. For the 11 remaining values the
condition $(*)$ holds if we take $\kappa=500$ and $C=10^{187}$; then the maximum upper bound 
$(**)$ is 231.

\emph{In this way we finally obtain the upper bound $H\leq 231=:K_1$, valid for all choices 
of $\al$ and $(i_0,j,k)$ mentioned in the beginning of this section.}
\\
\noindent
{\bf Computation time}. At this stage, the computation time was less than half of a minute.
\subsubsection{Further reduction and final stage of resolution}
             \label{subsubsec final}
We repeat the $p$-adic reduction process of Subsection \ref{subsubsec padic reduc} with $K_0\leftarrow K_1=231$ 
and $N_0\leftarrow N_1=207$. This affects $W$ and, consequently, $\kappa$ and $m$ in 
(\ref{eq choose m}), which now becomes ``very small''. Thus, we obtain the new bound $n_1\leq N_2=25$,
and this took less that 1 minute.
\\  \noindent
Next, applying the reduction process of Subsection \ref{subsubsec padic reduc} with $K_0\leftarrow K_1=231$ 
and $N_1\leftarrow N_2$, implies $H\leq K_2=41$; this took a few seconds.

A third $p$-adic reduction step can improve a little bit the upper bound for $n_1$. The process of 
Subsection \ref{subsubsec padic reduc} with $K_0\leftarrow K_2=41$ and $N_0\leftarrow N_2=25$ implies 
$n_1\leq N_3=21$, and this took around 1 second. 
Although we can make a further reduce to the bound of $H$, as well, and obtain $H\leq 34$, we will not
use this.

Actually, we prefer to solve a set of Thue equations (\ref{Thue eqs}), with right-hand side
$c\in\{-2^5 3^4 5^{z_1}11^{z_2} : 0\leq z_1\leq 27,\, z_2=0 \;\mbox{or}\; 3\leq z_2\leq 21\}$, 
using {\sc magma}'s implementation of Bilu \& Hanrot's method \cite{BiHa}. Remember that, as already 
mentioned a few lines below (\ref{Thue eqs}), no solutions exist when  
$c\in\{-2^5 3^4 5^{z_1}11^{z_2} : 0\leq z_1\leq 27,\, 1\leq z_2\leq 2\}$.
Thus, we are left with $28\times 20$ Thue equations, using the above implementation. This is the most 
expensive task; it took us about $7662\,\mbox{secs}\approx 2\mbox{h}\,7'\,42''$. 
No solutions were found, hence, we have the following result:
\begin{proposition}
      \label{prop n=5}
There are no solutions to the equation \emph{(\ref{eq ThueMahler5})}, hence, by 
Subsection \emph{\ref{subsec reduce to a TM eq}}, equation \emph{(\ref{eq title})} with $abx$ odd and 
$n=5$ has no solutions.  \noindent\proofend
\end{proposition} 
\noindent
Now, Propositions \ref{prop n<7} and \ref{prop n=5} complete the proof of Theorem \ref{thm main theorem}.
\\[2mm] \noindent
{\bf Remark}: Taking into account the computation cost at previous stages from Subsections 
\ref{subsubsec 1st ub} through \ref{subsubsec real reduc}, we see that the total computation time for the 
needs of Subsection \ref{subsec solution of T-M eq} is less than $2\mbox{h}.30'$.  
%%
%%%%%%%%%%%%%%%%% The Appendices %%%%%%%%%%
\begin{appendices}
\section{}
      \label{Append working padically}

\vspace*{-5mm}
{\large {\bf Working $p$-adically. Some general facts.}}
\\ \noindent
In this appendix we combine several facts which are scattered in the literature. Our basic references
are \cite{BoSh}, \cite{Ca}, \cite{FT}, \cite{Ko}, \cite{Mi}.

Let $p$ be a rational prime. For every non-zero $x\in\Q$ we denote by
$\ekth{p}{x}$ the exponent with which $p$ appears in the prime factorization
of $x$ and, as usually, the $p$-adic absolute value of $x$ is defined by
$|x|_p=p^{-\ekth{p}{x}}$. We set, by convention, $\ekth{p}{0}=-\infty$, so that
$|0|_p=0$. For $x\in\Q$, we also define $\ord{p}{x}=\ekth{p}{x}$.
\\ \noindent
This extends to $\Q_p$. 
If $x\in\Q_p$ and we write $x$ in the standard $p$-adic representation 
$x=\sum_{i=N}^{\infty} a_ip^i$ ($N\in\Z$, the $a_i$'s are integers with 
$0\leq a_i<p$ and $a_N\neq 0$), then we define $\ord{p}{x}=N$ and 
$|x|_p=p^{-N}$. Clearly, in the special case $x\in\Q$, these definition agree 
with those given above. 
\\ \noindent
More generally, if $x\in E_p$, where $E_p$ is a finite extension of $\Q_p$,
of degree, say $d$, and $t^d+b_{d-1}t^{d-1}\cdots+b_1t+b_0$ is the 
characteristic polynomial of $x$ with respect to the extension $E_p/\Q_p$, 
($b_i\in\Q_p$ for $i=0,\ldots,d-1$), then
\begin{equation} \label{eq order formula}
\ord{p}{x}=\frac{1}{d}\,\ord{p}{b_0}=\frac{1}{d}\,\ord{p}{\Nm{E_p}{\Q_p}(x)}
\quad\mbox{and}\quad |x|_p=p^{-\ord{p}{x}}.
\end{equation}
These definitions are independent from $E_p$;
in particular, they coincide with the definitions of $\ord{p}{x}$ and $|x|_p$ 
given at the beginning with $x\in\Q_p$. 

Now we adopt a different point of view. Let 
\[
E=\Q(\xi), \mbox{ where $g(\xi)=0$ and $g(t)\in \Q[t]$ is monic and irreducible.}
\]
We denote by $\O_E$ the maximal order of $E$. Let   
\begin{equation} \label{eq ideal factorization of <p>}
p\,\O_E=\idp_1^{e_1}\cdots\idp_m^{e_m}
\end{equation}
be the factorization of the principal ideal $p\,\O_E$ into prime ideals of $E$,
where the $\idp_i$'s above are distinct and ramification index $e_{E/\Q}(\idp_i)=e_i>0$ for every 
$i=1,\ldots,m$; we also denote by $f_i$ the residual degree $f_{E/\Q}(\idp_i)$.
\\ \noindent
For every $x\in E$ and every $\idp_i$ we denote by $\ekth{\idp_i}{x}$ the 
exponent of $\idp_i$ in the prime ideal factorization of $x\,\O_E$; 
in particular, $\ekth{\idp_i}{p}=e_i$. If $\idp_i$ is a principal ideal,
say $\idp_i=\pi\,\O_E$, then we write $\ekth{\pi}{x}$.

The polynomial $g(t)$ factorizes into $m$ distinct irreducible polynomials of
$\Q_p[t]$:
\begin{equation} \label{eq factor g(t) over Q_p}
g(t)=g_1(t)\cdots g_m(t).
\end{equation}
Let
\[
E_{\idp_i}=\Q_p(\xi_{\idp_i}), \;\; \mbox{where $\xi_{\idp_i}$ is defined by $g_i(\xi_{\idp_i})=0$;}
\]
actually, $E_{\idp_i}$ is the completion of $(E,|\cdot |_{\idp_i})$, where $|\cdot |_{\idp_i}$ is the 
absolute value of $E$ corresponding to the additive valuation $\ekth{\idp_i}{\cdot}$. 
There is a natural embedding $E\stackrel{\phi_i}{\hookrightarrow} E_{\idp_i}$ mapping 
$\xi$ to $\xi_{\idp_i}$, which allows us to view $E$ as a subfield of $E_{\idp_i}$.
The typical element $x(\xi)\in E$ (where $x[t]\in\Q[t]$) can be viewed as an element of  $E_{\idp_i}$
if we identify $x(\xi)$ with $\phi_i(x(\xi))$. Formally, this means that we view $x(\xi)$ as the element 
$x(\xi_{\idp_i})=x(t)+g_i(t)\Q_p[t] \in \Q_p[t]/g_i(t)\Q_p[t]$.
Then,  according to (\ref{eq order formula}),
\begin{equation} \label{eq def ordp}
\ord{p}{x(\xi)}=\frac{1}{[\Q_p(\xi_{\idp_i}):\Q_p]}\,\ord{p}{\Nm{\Q_p(\xi_{\idp_i})}{\Q_p}(x(\xi_{\idp_i}))}
\quad\mbox{and}\quad |x(\xi)|_p=p^{-\ord{p}{x(\xi)}}.
\end{equation} 
The above discussion makes clear that the value of $\ord{p}{x(\xi)}$ depends on $\idp_i$. 
Consequently, if $i\neq j$ and $x(\xi)\in E$, then, the value of $\ord{p}{x(\xi)}$ may vary, depending on 
whether we view $E$ as a subfield of $E_{\idp_i}$ or of $E_{\idp_j}$. 

The enumeration of the $\idp_i$'s in (\ref{eq ideal factorization of <p>})
and the $g_i$'s in (\ref{eq factor g(t) over Q_p}) can be done in 
such a way that 
\begin{equation}\label{eq correspond p_i with g_i}
\deg g_i = [E_{\idp_i}:\Q_p]=e_if_i\quad\mbox{and}\quad
\ekth{\idp_i}{x(\xi)}=e_i\cdot\ord{p}{x(\xi)}.
\end{equation}
The second relation above implies that, for the typical element $x(\xi)\in E$ (where 
$x[t]\in\Q[t]$), the following is true: $x(\xi)$ is divisible by $\idp_i$ iff $\ordp(x(\xi_{\idp_i}))>0$ which,
in turn, is equivalent to the statement that the constant term of the characteristic polynomial of 
$x(\xi_{\idp_i})$ over $\Q_p$ has positive $\ordp$.
Having established this enumeration, we have for every $x(\xi)\in E$:
\begin{equation} \label{eq nu_pi(x)}
\ekth{\idp_i}{x(\xi)}=\frac{e_{E/\Q}(\idp_i)}{[E_{\idp_i}:\Q_p]}\,\ord{p}{\Nm{E_{\idp_i}}{\Q_p}(x(\xi_{\idp_i}))}
    = e_{E/\Q}(\idp_i)\cdot\ord{p}{x(\xi)}.
\end{equation}

In practice, the above mentioned correspondence $\idp_i\leftrightarrow g_i$ is carried out by a 
{\sc magma} routine \label{page my magma routine} which we wrote based on the following: 
For $j=1,\ldots,m$, consider the ``two-element representation'' of $\idp_j$, namely, 
$\idp_j=p\,\O_E +h_j(\xi)\,\O_E$, where $h_j(t)\in\Q[t]$, and fix any $i\in\{1,\ldots,m\}$. 
For $j=1,\ldots,m$, compute $\phi_i(h_j(\xi))=h_j(\xi_{\idp_i})$ and the characteristic polynomial $\chi_j(t)$ of 
$h_j(\xi_{\idp_i})$ with respect to the extension $E_{\idp_i}/\Q_p$. 
For exactly one index $j$ the $\ordp$ of the constant term of $\chi_j(t)$ is 
positive. The polynomial $g_j(t)$, for this specific $j$, corresponds to the ideal $\idp_i$. This we do for any 
$i=1,\ldots,m$ and we establish the one-to-one correspondence 
$\{\idp_1,\ldots,\idp_m\}\leftrightarrow \{g_1,\ldots,g_m\}$. By permuting the
indices of $g_1,\ldots,g_m$, if necessary, we establish the one-to-one correspondence 
$\idp_i\leftrightarrow g_i$ which satisfies (\ref{eq correspond p_i with g_i}) and (\ref{eq nu_pi(x)}).
\section{}
             \label{Append field K}

 \vspace*{-5mm}            
{\large {\bf Working in $K$}.}
 \\[5pt] \noindent             
In order to apply the method of \cite{TdW} we need to work in an extension $K$ of $F$ in which $g(t)$ 
has  at least three distinct roots. Thus, in general, we do not need the whole splitting field of $g(t)$ 
over $F$. In our case, however, the Galois group of $g(t)$ is of order $20$, which implies that 
$K$ is the splitting field of $g(t)$ over $F$. Using {\sc magma} we find out that $K=\Q(\om)$, where
$\om$ is a root of the polynomial
{\small
\begin{equation*}
\begin{split}
G(t) & = \\
 & t^{20} + 780t^{19} + 248030t^{18} + 39929580t^{17} + 3046440525t^{16} + 
                                                                   18210793968t^{15}- 13729990391320t^{14}   \\
& 
+ 752551541981520t^{13} + 8605950990819730t^{12} + 1708764818389209000t^{11}  \\
&  
+23308084571944423284t^{10} - 1404817549102176551640t^9 -35442768652652017430190t^8 \\
& 
+ 375805034836819117590960t^7+16191084883780784798260200t^6 \\
&
+ 30210122048192693893581552t^5 -2113554835538935196795743635t^4 \\ 
&
- 12364486598313473303834175060t^3 + 25061666765667764525027943390t^2 \\
& + 278757784774895111708136641100t +  427756623168133431059207412321. 
\end{split}
\end{equation*}
}
Of course, $K$ is a Galois extension.
By $\O_K$ we denote the maximal order of $K$ and by $\O_F$ the maximal order of $F$.
The roots  $\cnj{\th}{i}$ ($i=1,\ldots,5$) of $g(t)$ are polynomial expressions of $\om$ 
with rational coefficients. 
Thus $\cnj{\th}{i}=\th_i(\om)\in\Q[\om]$, where $\th_i(t)\in\Q[t]$ ($i=1,\ldots,5$). 
Each $\cnj{\th}{i}$ corresponds to an embedding  $\psi_i: F\hookrightarrow K$  
 characterized by $\psi_i(\th)=\th_i(\om)$. Then we can view $F$ as subfield of $K$ in five ways, by
identifying $F$ with $\psi_i(F)$. For our computations we can arbitrarily choose the embedding $\psi_i$,
but, once we choose it, we must keep it fixed.  {\sc Magma} computes the embeddings $\psi_i$; 
 rather arbitrarily, when it considers $F$ as a subfield of $K$, our {\sc magma} session implicitly uses the 
 embedding $\psi_5$;   
 i.e. as a subfield of $K$, $F$ is identified with $\psi_5(F)$. For simplicity, we set $\psi=\psi_5$  
and avoid the use of superscript/subscript indicating  the conjugation. Then, extending a prime ideal $\idp$ 
of $\O_F$ to the ideal $\idp\,\O_K$ means, the following: Let $\idp=p\,\O_F+h(\th)\,\O_F$, where 
$h(\th)\in\Q[\th]$. Then $\idp\,\O_K= p\,\O_K+h(\psi(\th))\,\O_K$.  {\sc Magma} computes  
$\psi(\th)=\psi_5(\th)=\cnj{\th}{5}=\th_5(\om)$:
{\small
%\begin{equation*}
\begin{align*}
\th_5(\om)  = & (109949833761153867182006233162830218735318443\,\om^{19} \\
  &
             +85209181362831884900908863253209941719229099772\,\om^{18}  \\
                  &
 +26843948030443532446996152295591813313759057989280\,\om^{17}   \\
 &
 +4256032326064110742980882546282193567817773030947507\,\om^{16}   \\
&  
+313740861737738480045017056015360916332191517809963652\,\om^{15}  \\
&
 +448056667342462133472342261671551993343950346279488008\,\om^{14}   \\
 & 
-1510732666058952542367528947242709391406430475149001889048\,\om^{13}  \\
&
-75203582469112061453856831841769761124516215292521423604700\,\om^{12}  \\
&
+1315665925830026167841281970821852638550108379556176872190482\,\om^{11} \\
& 
+181102009414236143682774128780867507407780245434648419322814672\,\om^{10}  \\
& 
+1669796300877028395254905627635289010395439563221417145651022616\,\om^9 \\
&   
-162206235802302893664200291374667388423732969287959076428108697406\,\om^8  \\
&\mbox{\small \sl (continued on next page)}
\end{align*}
\begin{align*}
&\mbox{\small \sl (continued from previous page)} \\
 &
-3096949958783354342722100732857399558416026343530279487756804602460\,\om^7  \\
&
  +56131367252940396214847119246725578907696249119733671305047923800344\,\om^6    \\
&
 +1504180072707173996628600619493969793943259186258161149867325499759064\,\om^5  \\
&
   -3870840446117907139295078483933668999054374262639172467797695558733756\,\om^4  \\
&
  -213744087631757677190145693131224085047423871295167330938194597435156069\,\om^3   \\
&
  -352243194516613905751693267069064334504141157799975445500242914720886380\,\om^2   \\
&
   +4494640381465135578890041723352987174508368793027744613522589459354871816\,\om   \\
&
 +10018536620467342924854600525808560481514235129562992768169198074640026531)/\\
& \hspace*{15mm}
             23750735994552570259738911035979918362684670692062523272986624000000000.
\end{align*}
%\end{equation*}
}
For the ideal $\idp=\pi_{113}\,\O_F$, obviously, we can take $h(\th)=\pi_{113}$ 
(see page \pageref{page def pi_113}). Then 
{\small
\begin{align*}
\psi(h(\th))=\psi(\pi_{113})  = & 
 (4378482585825381431566239057028627160970028037\,\om^{19} \\
 &
+3396229499207087892836807551803460134822567613017\,\om^{18} \\
&
+1071214373586764983689778732938427724586732483076327\,\om^{17} \\
&
+170145388685753978342778727714777154218418254512908043\,\om^{16} \\
&
+12588645126248857287103560864595388726070407650893086772\,\om^{15}  \\
&
+23291648519661572328486403102949754117988678982445451332\,\om^{14} \\
&
-60338282775505500617169124951383242755661107414046130424148\,\om^{13} \\
& 
-3031547893595987199227070326893238313728050147362743168899908\,\om^{12} \\
&
+51530222400182570158994336776484627607851706804827046320630198\,\om^{11} \\
&
+7282924550357450053556161972634939332060481531424832236222347982\,\om^{10} \\
&
+69407725895185343408814814665447944190629387800816493362801270434\,\om^9 \\
&
-6522081891215060359860869744209143119462208244901441175704137508646\,\om^8 \\
&
-126498537591831789061901120743476135204565837943425968840309704817756\,\om^7 \\
&
+2266682290068533497159436879932180463046527604540137185016229550764276\,\om^6 \\
&
+61296290267956300004998582874965020844943081524921068036926148986534364\,\om^5 \\
&
-165499773355462341124868225890607890561889109853343804837744925708526164\,\om^4 \\
&
-8698248530644221957092224575657663153179204605475867759713929659343969939\,\om^3 \\
&
-11520188167294592917194821184972443593645584955137753402305064016274508639\,\om^2 \\
&
+173158801650947255768584220314522267452080115712427401970346621418344586479\,\om \\
&
+378768220799897336478183837425696328863796308850264844412883432915449692851)/ \\
&
\hspace*{15mm}1520047103651364496623290306302714775211818924292001489471143936000000000.
\end{align*}
}
Our {\sc magma} routine, mentioned below the equation (\ref{eq nu_pi(x)}), returns the
factorization $11\O_K=\prod_{i=1}^{10}\idP_i^2$, where the residual degree $f_{K/\Q}(\idP_i)=1$ for every 
$i=1,\ldots,10$.
Moreover, for every $i=1,\ldots,10$, the routine computes:
\\ \noindent
\textbullet\, An element $h_i\in\O_K$, such that $\idP_i=11\O_K+h_i\O_K$. 
\\ \noindent
\textbullet\,
The factorization of $G(t)=\prod_{i=1}^{10}G_i(t)$ into irreducible polynomials over $\Q_p$. 
For $i=1,\ldots,10$, the irreducible polynomial $G_i(t)\in\Q_p[t]$ corresponds to $\idP_i$ in the sense 
explained in Appendix \ref{Append working padically}. As expected, $\deg G_i(t)=2$ for every 
$i=1,\ldots,10$.
\\ \noindent
In table \ref{Table factorization in K} we give the data mentioned in the above two ``bullets''. 
The element $h_i$ is identified with a $20$-tuple: $h_i=(c_{i1},\ldots,c_{i20})$ means that 
\[
h_i=\sum_{j=1}^{20} c_{ij}\be_j,
\]
where $\be_1,\ldots,\be_{20}$ is an integral basis of $K/\Q$, explicitly calculated by {\sc magma}.
For the polynomial $G_i(t)$ we write $G_i(t)=(\ga_{i1},\ga_{i0})$, by which we mean that 
$G_i(t)=t^2+\ga_{i1}t+\ga_{i0}$.  In the columns of the $\ga_{ij}$'s we write their $11$-adic approximations,
(rational integers) with precision $O(11^{10})$.
\begin{center}
\tablecaption{$\quad 11\O_K=\prod_{i=1}^{10}\idP_i^2$ }   
                                                      \label{Table factorization in K}  
\tablefirsthead{ \\[-10mm] \hline
\multicolumn{1}{|c|}{\tbsp } & 
\multicolumn{1}{|c|}{$\idP_i=11\O_K+h_i\O_K$} &
\multicolumn{2}{|c|}{$K_{\idP_i}=\Q_{11}[t]/\langle t^2+\ga_{i1}t+\ga_{i0}\rangle$} 
\\ \hline
   $i$ & $h_i$  & $\ga_{i1}$ & $\ga_{i0}$ \\ \hline\hline
               }
\tablehead{ 
\multicolumn{2}{c}{\multicolumn{4}{|l|}{{\small \sl continued from previous page}}} \\
\hline
\multicolumn{1}{|c|}{\tbsp } & 
\multicolumn{1}{|c|}{$\idP_i=11\O_K+h_i\O_K$} &
\multicolumn{2}{|c|}{blabla} 
\\ \hline
   $i$ & $h_i$  & $\ga_{i1}$ & $\ga_{i0}$ \\ \hline\hline
          }
\tabletail{\hline \multicolumn{2}{|r|}{\small\sl continued on next page}\\ \hline}
\tablelasttail{ }
 \begin{supertabular}{|c|c|c|c|}
$1$ & $(5, 8, 10, 5, 3, 3, 7, 0, 5, 5, 0, 3, 2, 10, 0, 2, 10, 2, 1, 7)$ 
& $11244468595$ & $6668815422$ 
 \\ \hline
$2 $ & $(6, 2, 0, 0, 7, 10, 4, 0, 6, 6, 4, 5, 2, 10, 0, 2, 10, 2, 1, 7)$ 
& $ -169320583$ & $10491152974$  
 \\ \hline
$ 3$ & $(6, 9, 8, 0, 6, 6, 1, 5, 2, 0, 6, 5, 2, 10, 0, 2, 10, 2, 1, 7)$ 
&$ -9236124994$ & $ 5583083423$  
 \\ \hline
$ 4$ & $(1, 5, 8, 3, 7, 5, 10, 6, 9, 0, 9, 3, 2, 10, 0, 2, 10, 2, 1, 7)$ 
&$ -6126278749$ & $ 7582171800$  
 \\ \hline
$ 5$ & $(2, 1, 10, 4, 6, 7, 8, 6, 2, 9, 1, 3, 7, 6, 5, 10, 7, 2, 10, 7)$ 
&$10744341441 $ & $-10666285673 $  
 \\ \hline
$ 6$ & $(8, 4, 10, 3, 5, 4, 3, 4, 4, 1, 5, 3, 4, 0, 8, 2, 1, 0, 8, 7)$ 
&$ -3779293982$ & $290904043 $  
 \\ \hline
$ 7$ & $(10, 8, 2, 1, 0, 9, 10, 4, 8, 4, 4, 8, 10, 2, 1, 7, 9, 4, 8, 7)$ 
&$-669447737 $ & $7802303026 $  
 \\ \hline
$8 $ & $(10, 1, 10, 2, 6, 10, 6, 7, 2, 1, 1, 0, 5, 1, 6, 10, 5, 3, 8, 7)$ 
&$ -12083236915$ & $-10106323842 $  
 \\ \hline
$ 9$ & $(8, 7, 9, 1, 9, 10, 10, 3, 3, 0, 8, 10, 10, 4, 6, 6, 6, 8, 8, 7)$ 
&$10744341441 $ & $9625552201 $  
 \\ \hline
$10 $ & $(6, 0, 10, 3, 6, 2, 7, 9, 3, 7, 6, 1, 8, 10, 0, 1, 9, 1, 2, 9)$ 
&$-669447737 $ & $-12489534848 $  
\\ \hline\hline
 \end{supertabular}
  \end{center}

Now, for $\idp=\pi_{113}\O_F$ we have to know the factorization of $\idp\,\O_K$; of course, the prime 
ideals of $\O_K$ in this factorization belong to $\{\idP_1,\ldots,\idP_{10}\}$.  For this purpose it suffices
to compute $\ekth{\idP_i}{\psi(\pi_{113})}$ for $i=1,\ldots,10$. This we do easily using {\sc magma}.
We find out that $\ekth{\idP_i}{\psi(\pi_{113})}=0$ for $i=1,\ldots,8$ and $\ekth{\idP_i}{\psi(\pi_{113})}=2$
for $i=9,10$; hence $\pi_{113}\,\O_F=\idP_9^2\,\idP_{10}^2$.  According to the above we put
$\idP=\idP_9$, so that
\label{page choice of mathfrakP}
\begin{equation}  \label{eq e,f of idP}
e_{K/\Q}(\idP)=2, \quad f_{K/\Q}(\idP)=1
\end{equation}
\[
G_{\idP}(t)=G_9(t)=t^2+\ga_{91}t+\ga_{90}. 
\]
Working $p$-adically in $K$ means working in 
$K_{\idP}=\Q_{11}(\om_{\idP})\cong\Q_{11}[t]/\langle G_{\idP}(t)\rangle $, where each root $\cnj{\th}{i}$ 
($i=1,\ldots,5$) is identified with $\th_i(\om_{\idP})$. 
\\
We have $g_3(\th_5(\om_{\idP}))=0$. Indeed, we have the following commutative diagram of 
monomorphisms:
\[
\begin{CD}
F     @>{\psi=\psi_5}>> K   \\[-5pt]
@VV\phi V               @VV\Phi V  \\[-5pt]
F_{\idp} @>{\Psi}>> K_{\idP}
\end{CD}
\]
Here $F_{\idp}=\Q_p(\th_{\idp})$, where 
$\th_{\idp}=7 + 2\cdot 11 + 2\cdot 11^2 + 10\cdot 11^3 + 7\cdot 11^4 +\cdots \in\Q_{11}$, is the root of 
$g_3(t)$.\footnote{We remind that $g_1(t), g_2(t), g_3(t)$ are defined at the beginning of 
Subsection \ref{subsubsec working 11adically}. } 
The natural embeddings $\phi$ and $\Phi$ are in accordance with the general discussion a few lines below
the relation (\ref{eq factor g(t) over Q_p}). Thus, $\phi(\th)=\th_{\idp}$, $\Phi(\om)=\om_{\idP}$ and,
consequently,
\[
\Psi(\th_{\idp})=\Phi\circ\psi_5\circ\phi^{-1}(\th_{\idp})=\Phi\circ\psi_5(\th) =\Phi(\th_5(\om))=
\th_5(\Phi(\om))=\th_5(\om_{\idP}).
\]
Therefore, $g_3(\th_5(\om_{\idP}))=g_3(\Psi(\th_{\idp}))=\Psi(g_3(\th_{\idp}))=\Psi(0)=0$.
Further, using {\sc magma} we see that the ($11$-adic) roots of $g_1(t)$ are $\th_i(\om_{\idP})$ with 
$i=2,4$, and the roots of $g_2(t)$ are $\th_i(\om_{\idP})$ with $i=1,3$.

In Subsection \ref{subsubsec n1<c13(logH+c14)},  where we view (\ref{eq (12) of TdW})  as a relation in 
$K_{\idP}$ (which simply means that we apply $\Phi$ to (\ref{eq (12) of TdW})) we will choose $i_0=5$, 
$j=1$ and $k=3$, following the instructions at bottom of p.~235 of \cite{TdW}. 
\end{appendices}

\end{document}